\documentclass{elsarticle}[preprint,review,10pt]

\journal{Journal of Computational Physics}

\usepackage{tikz}
\usetikzlibrary{arrows.meta, positioning}
\usepackage{graphicx} 

\usepackage{graphicx} 
\usepackage{algorithm}
\usepackage{algorithmic}
\usepackage{booktabs} 
\usepackage{siunitx} 
\usepackage{caption} 
\usepackage{multirow} 
\usepackage{threeparttable}

\usepackage[utf8]{inputenc}
\usepackage{graphicx}
\usepackage{bm}
\usepackage{comment}
\usepackage{psfrag}
\usepackage{latexsym,amsmath,amsfonts,amscd,amsthm}
\usepackage{changebar}
\usepackage{color}
\usepackage{bm}
\usepackage{tikz}
\usepackage{multirow}
\usepackage{xcolor}
\usepackage{hyperref}
\usepackage{setspace}
\usepackage{amssymb}
\usepackage{mathrsfs}
\usepackage{caption}
\usepackage{subcaption}
\usepackage[title]{appendix}
\usepackage{ulem}
\usepackage{tikz}
\usetikzlibrary{arrows,backgrounds}

\newlength{\spse}
\newcommand{\Dt}{\Delta t}

\newcommand{\BB}{\mathbf{B}}

\newcommand{\BU}{\mathbf{U}}
\newcommand{\BV}{\mathbf{V}}

\newcommand{\BT}{\mathbf{T}}

\newcommand{\BI}{\mathbf{I}}

\newcommand{\BR}{\mathbf{R}}

\newcommand{\WBA}{\widetilde{\mathbf{A}}}
\newcommand{\WBC}{\widetilde{\mathbf{C}}}
\newcommand{\bG}{\mathbf{G}}

\newcommand{\bv}{\boldsymbol{\upsilon}}

\newcommand{\bx}{\mathbf{x}}

\newcommand{\bc}{\mathbf{c}}
\newcommand{\bb}{\mathbf{b}}

\newcommand{\wbb}{\widetilde{\mathbf{b}}}

\newcommand{\bn}{\boldsymbol{n}}

\newcommand{\half}{\frac{1}{2}}
\newcommand{\brho}{\boldsymbol{\rho}}
\newcommand{\bg}{\boldsymbol{g}}

\newcommand{\bq}{\boldsymbol{q}}

\newcommand{\wbG}{\widetilde{\boldsymbol{G}}}

\newcommand{\bff}{\boldsymbol{f}}

\newcommand{\BC}{{\bf{C}} }

\newcommand{\BM}{{\bf{M}} }
\newcommand{\BH}{{{\bf{H}}}}
\newcommand{\BD}{{\bf{D}}}

\newcommand{\BS}{\mathbf{S}}
\newcommand{\BSigma}{{\mathbf{\Sigma}}}

\newcommand{\by}{{\bf{y}}}

\newcommand{\rT}{{\rm{T}}}

\newcommand{\pzc}{\textcolor{blue}}

\title{On-the-Fly ROM-Based Acceleration of SI-DSA for Implicit Time Marching of the Radiative Transfer Equation}

\makeatletter \def\ps@pprintTitle{} \makeatother

\begin{document}

\author{
Ningxin Liu, Zhichao Peng
}

\ead{nxliu@ust.hk,pengzhic@ust.hk}

\affiliation{organization={Department of Mathematics, The Hong Kong University of Science and Technology},
            addressline={Clear Water Bay, Kowloon},
            city={Hong Kong},
            country={China}
}

\begin{abstract}

In implicit time marching of the radiative transfer equation (RTE), the resulting linear systems are commonly solved using source iteration with diffusion synthetic acceleration (SI-DSA). Despite its widespread success, the performance of the DSA preconditioner may deteriorate when the RTE cannot be well approximated by its diffusion limit. Moreover, classical SI-DSA does not exploit low-rank structures of the solution manifold across time steps when the solution evolves smoothly.

To address these limitations, we develop an on-the-fly reduced-order-model (ROM)-based acceleration for SI-DSA in implicit time marching of the RTE. Instead of relying on a diffusion approximation, the proposed approach constructs ROMs directly from the underlying kinetic formulation while exploiting low-rank structures in the temporal evolution of the solution. The method is fully offline-free and constructs ROMs to enhance both initial guesses and preconditioners on the fly during time marching. To handle streaming solution data, we design efficient and memory-lean ROM construction and adaptive update strategies based on dynamical mode decomposition, incremental low-rank singular value decomposition, and error indicators. Numerical experiments demonstrate that the proposed method consistently accelerates implicit time marching, delivering 
$1.4\times$ to $2.0\times$
speedup over classical SI-DSA while incurring only marginal overhead for ROM construction and updates.
\end{abstract}


\maketitle

\section{Introduction}

The radiative transfer equation (RTE) is a kinetic equation that models particles such as neutrons or photons propagating through and interacting with a background medium. It provides basic models in applications including nuclear engineering, medical imaging, astrophysics and atmospheric science. This equation is featured by its multiscale property: the behavior of particle system governed by it is transport dominant in weak scattering regions and becomes diffusion dominant in strong scattering regions. To efficiently capture this multiscale behavior, asymptotic preserving implicit time marching is widely applied \cite{larsen2009advances}.   
An efficient iterative solver for the resulting linear system is crucial to the overall efficiency.

Source iteration with diffusion synthetic acceleration (SI-DSA) is a widely applied iterative solver for the linear systems arising in implicit time marching of the RTE \cite{adams2002fast,larsen2009advances}. DSA \cite{kopp1963synthetic,adams1992diffusion,wareing1993new} acts as a crucial preconditioning step that accelerates the convergence of source iteration. Its core idea is to introduce a correction to the iterative solution by approximating an ideal kinetic correction equation with its diffusion limit. 

Despite the great success of classical SI-DSA, two limitations remain. First, the effectiveness of DSA relies on the validity of the diffusion approximation and may deteriorate in problems far from the diffusion regime \cite{ren2019fast}. Second, when the solution varies smoothly in time, the solution manifold across time steps often exhibits low-rank structures, which are not exploited by classical SI-DSA.

Data-driven reduced order models (ROMs) have been developed for RTE in \cite{buchan2015pod,tencer2017accelerated,choi2021space,peng2024micro,tano2021affine,coale2023reduced,matsuda2025reduced} and can be applied to address these two limitations. Rather than relying on the diffusion approximation, ROMs can be constructed directly from the original kinetic formulation by leveraging low-rank structures of the solution manifold across time extracted from solution data.

In our previous work \cite{peng2024reduced,peng2025flexible,tang2025synthetic}, a family of ROM-based preconditioners was developed for steady-state parametric RTE. However, extending these approaches to implicit time marching presents new challenges. Unlike parametric problems where training solution data are available in advance, solution data are computed on the fly during the implicit time marching. In addition, the ROMs in \cite{peng2024reduced,peng2025flexible,tang2025synthetic} are constructed from data of the microscopic particle distribution, leading to high memory cost. 

To overcome these  challenges, we design a three-phase on-the-fly ROM-based acceleration for implicit time marching based on a memory-efficient reformulation of the RTE, whose unknown is the macroscopic density rather than the microscopic distribution.

\begin{enumerate}
\item \textbf{Phase 1, construct a ROM to improve initial guess.} We incrementally construct a ROM for the memory-efficient formulation of the RTE to provide an improved initial guess using dynamical mode decomposition (DMD) \cite{schmid2010dynamic} together with incremental low-rank singular value decomposition (SVD) \cite{brand2006fast}.

\item \textbf{Phase 2, construct a ROM to enhance preconditioner.} We use the ROM-enhanced initial guess for SI-DSA and incrementally build a second ROM for the kinetic correction equation to enhance the preconditioner. In addition, the ROM for the initial guess is adaptively updated based on an error indicator to maintain its accuracy.

\item \textbf{Phase 3, full ROM-enhanced SI-SA.}  We apply both the ROM-enhanced initial guess and a hybrid preconditioner similar to \cite{peng2024reduced,tang2025synthetic}, which exploits ROM to reduce low-frequency error components and DSA to eliminate high-frequency error components. Moreover, these ROMs are adaptively updated based on error indicators to maintain their accuracy.
\end{enumerate}
The key components of the proposed method are summarized as follows.
\begin{enumerate}
\item \textbf{On-the-fly ROM construction and adaptation for streaming solution data.}
To accommodate the streaming solution data generated during implicit time marching, we construct and adaptively update ROMs using dynamical mode decomposition (DMD) \cite{schmid2010dynamic} together with incremental low-rank SVD \cite{brand2006fast}, guided by error indicators and spectral information. These techniques enable efficient ROM construction and adaptation without requiring precomputed training data.

\item \textbf{Efficient and memory-lean ROM construction for the macroscopic formulation of the RTE.}
Although the projection-based ROM in \cite{behne2022minimally} also utilizes the memory-efficient formulation for the macroscopic density, it constructs the reduced system by projecting a sweep-based streaming operator with computational cost $\mathcal{O}(N_{\bx}N_{\bv}r^2)$, where $N_{\bx}$,  $N_{\bv}$, and $r$ denote the degrees of freedom in the physical space, the angular space, and the rank, respectively. In contrast, DMD avoids this expensive projection and constructs the reduced system in a purely data-driven manner with cost $\mathcal{O}(N_{\mathbf{\bx}}r^2)$. 

\item \textbf{Sequential construction of ROMs for the initial guess and the preconditioner.}
As pointed out in \cite{tang2025synthetic}, sequential construction of ROMs for the initial guess and the preconditioner is crucial to the effectiveness of ROM-enhanced preconditioners, since the error equation for SI-DSA implicitly depends on the initial guess. In this papere, we extend the sequential ROM construction for steady state parametric RTE in \cite{tang2025synthetic} to the streaming solution data generated during time marching.
\end{enumerate}
Numerical tests demonstrate that the proposed on-the-fly ROM acceleration to SI-SA is able to achieve $1.4\times$ to $2.0\times$ speedup over SI-DSA by paying marginal overhead ($1\%$ or $2\%$) to construct and update ROMs.

To better contextualize our work, we note that other ROM-based accelerations have also been developed for iterative solvers of the RTE. A DMD-based preconditioner was proposed for the implicit time marching of thermal radiation in \cite{mcclarren2019calculating} as an alternative to DSA in situations where DSA is not applicable, for example when nonlinear positivity-fixing are applied. The key difference between \cite{mcclarren2019calculating} and our method is that the approach in \cite{mcclarren2019calculating} exploits low-rank structures across  iterations within the same time step, whereas our method leverages low-rank structures across different time steps. In fact, the two approaches are complementary and can be combined by replacing DSA in our hybrid preconditioner with the method in \cite{mcclarren2019calculating}. Beyond the SI framework, another ROM-based acceleration for RTE is proposed in  \cite{chen2021low} by building a fast boundary-to-boundary mapping for a domain decomposition solver via randomized SVD.

The paper is structured as follows. Sec.  \ref{sec:preliminary} provides a brief background review of numerical methods for time-dependent RTE  and  SI-DSA. Sec.  \ref{sec:ROM_construct} introduces our on-the-fly ROM-based acceleration of SI-DSA in implicit time marching.  Sec.  \ref{sec:numerical} presents a series of numerical experiments to evaluate the performance of the proposed approach. Sec.  \ref{sec:conclusion} draws  conclusions and outlines potential future work.

\section{Preliminary}
\label{sec:preliminary}
In this paper, we consider the following time dependent linear radiative transfer equation with one energy group, isotropic scattering,  and isotropic inflow boundary conditions:
\begin{subequations}
\label{eq:rte}
    \begin{align}   &\partial_{t} f(\bx,\bv,t)+\bv\cdot\nabla_{\bx}f(\bx,\bv,t) + \sigma_{t}(\mathcal{\mathbf{x}})f(\mathcal{\bx},\bv,t) = \sigma_{s}(\bx)\rho(\bx,t) + G(\bx),\;\; \bx\in \Omega_{\bx}, \;\; \bv\in\mathbb{S}^2,\;\; t\in[0,T],
    \label{eq:kinetic_1}\\
    &\rho(\bx,t) = \frac{1}{4\pi}\int_{\bv\in\mathbb{S}^2} f(\bx,\bv,t) d\bv,  \label{eq:governing}
    \\
    &f(\bx, \bv,0) = f_0,
    \quad f(\bx, \bv,t) = g(\bx), \;\; \bx\in\partial\Omega_{\bx},\;\; \bv\cdot \mathbf{n}(\bx)<0.
    \label{eq:boundary}
    \end{align}
    \label{eq:kinetic}
\end{subequations}
Here $f(\bx, \bv,t)$ is the particle distribution function of time $t$, position $\bx$, and angular direction $\bv$, $\rho(\bx,t)$ is the macroscopic density, $\sigma_s\geq0$ is the scattering cross section, $\sigma_t(\bx)\geq\sigma_s(\bx)\geq0$  is the total cross section, $\sigma_a(\bx)=\sigma_t(\bx)-\sigma_s(\bx)$ is the absorption cross section, $G(\bx)$ is an isotropic  source term, and $\mathbf{n}(\bx)$ is the outward normal direction of $\Omega_{\bx}$ at $\bx\in\partial\Omega_{\bx}$.

When the scattering cross section $\sigma_s\rightarrow\infty$,  the particle distribution converges to the macroscopic density, i.e. $f\rightarrow\rho$, governed by the diffusion limit:
\begin{equation}
\partial_{t}\rho(\bx,t)-\nabla_{\bx}\cdot\Big(\BD_{\bv}\nabla_{\bx}\rho(\bx,t)\Big)+\sigma_a(\bx)\rho(\bx,t)=G(\bx),\;\; \bx\in \Omega_{\bx},\;\; t\in[0,T],
\end{equation}
where
$
\mathbf{D}_{\bv} = \frac{1}{4\pi\sigma_s}\textrm{diag}\left(\int_{\mathbb{S}^2} \bv_x^{2}d\bv, \;\int_{\mathbb{S}^2} \bv_y^{2}d\bv, \;\int_{\mathbb{S}^2} \bv_z^{2}d\bv\right)=\frac{1}{3\sigma_s}I_{3\times 3}.$

In this section, we  first introduce discretizations for RTE in Sec. \ref{sec:discretization}, then present the resulting linear system in Sec. \ref{sec:mat-vec} and discuss Source Iteration with Synthetic Acceleration (SI-SA) to solve the resulting discrete system in Sec. \ref{sec:sisa}. Though focusing on $S_N$ upwind DG method, the proposed on-the-fly ROM acceleration techniques can potentially be applied to other discretizations. 

\subsection{Implicit discrete ordinates upwind discontinuous Galerkin method\label{sec:discretization}} 
Here, we outline our discretization for equation \eqref{eq:rte}.

\textbf{Discrete ordinates ($S_N$) angular discretization.}
In angular space, we apply the discrete ordinates ($S_N$) method \cite{pomraning2005equations}. Specifically, $S_N$ method solves RTE at a set of quadrature points $\{\bv_j\}_{j=1}^{N_v}\subset\mathbb{S}^2$, which corresponds to quadrature weights $\{\omega_j\}_{j=1}^{N_v}$ satisfying $\sum_{j=1}^{N_v}\omega_j=1$. The normalized integral term $\rho(\bx,t)$ in \eqref{eq:kinetic} can be approximated as:
\begin{equation*}
    \rho(\bx,t)  = \frac{1}{4\pi}\int_{\mathbb{S}^2} f(\bx, \bv,t) d\bv \approx \sum_{j=1}^{N_{\bv}} \omega_{j} f(\bx, \bv_j,t).
\end{equation*}
Then, $S_N$ discretization for  RTE (\ref{eq:kinetic}) is defined as:
\begin{subequations}
\label{eq:rte discrete}
    \begin{align}
    &(\partial_{t}+\bv_{j}\cdot\nabla_{\bx}+\sigma_t(\bx))f(\bx,\bv_j,t) = \sigma_{s}(\bx)\rho(\bx,t) + G(\bx),\;\; \rho(\bx,t) = \sum_{j=1}^{N_{\bv}} \omega_{j} f(\bx,\bv_{j},t), \label{eq:governing_discrete}
    \\
    &f(\bx, \bv_{j},t) = g(\bx), \;\;\bx\in\partial\Omega_{\bx},\;\; \bv_j\cdot \mathbf{n}(\bx)<0.
    \label{eq:boundary_discrete}
    \end{align}
\end{subequations}
In this paper, we apply Gauss-Legendre quadrature points for 1D slab geometry and Chebyshev-Legendre (CL) rule in higher dimensions. For brevity, details of the CL quadrature are deferred to \ref{appx:cl}. 

\textbf{Time discretization.}
As for time discretization, we consider the backward Euler method. Let $\Delta t$ denote the time step size. The discrete time levels are then given by $t^n:=n\Delta t$ for $n\in \mathbb{N}$.  The solution of $f(\bx, \bv_{j},t)$ at time $t^n$ is approximated by $f(\bx, \bv_{j},t^n)$. Then discretization scheme of Eq. (\ref{eq:kinetic_1}) is given by
\begin{equation*}
\frac{f(\bx,\bv_j,t^n)-f(\bx,\bv_j,t^{n-1})}{\Delta t}+\bv_{j}\cdot\nabla_{\bx}f(\bx,\bv_j,t^n)+\sigma_t(\bx)f(\bx,\bv_j,t^n) = \sigma_{s}(\bx)\rho(\bx,t^n) + G(\bx),
\end{equation*}
where $\rho(\bx,t^n) = \sum_{j=1}^{N_{\bv}} \omega_{j} f(\bx,\bv_{j},t^n)$. Rearrange this scheme as 
\begin{equation}
\label{eq:time_discrete}
\Big(\frac{1}{\Dt}+\bv_{j}\cdot\nabla_{\bx}+\sigma_t(\bx)\Big)f(\bx,\bv_j,t^n)=\sigma_{s}(\bx)\rho(\bx,t^n)+\frac{1}{\Dt}f(\bx,\bv_j,t^{n-1})+G(\bx).
\end{equation}
Though focusing on backward Euler method in this paper, the proposed acceleration techniques can be coupled with other implicit time marching scheme such as backward differentiation formula (BDF).

\textbf{Upwind discontinuous Galerkin spatial discretization.}
To capture the diffusion limit of RTE without resolving small mean free path, we apply high order upwind discontinuous Galerkin (DG) spatial discretization which is proved to be asymptotic preserving \cite{adams2001discontinuous,guermond2010asymptotic}. We consider 2D X-Y geometry with a rectangular computational domain $\Omega_{\bx}=[x_l,x_r]\times[y_l,y_r]$ and a partition of it with rectangular meshes $\mathcal{T}_{h} = \{\mathcal{T}_{i},\forall \mathcal{T}_{i} {\rm \;being\;rectangle}\}_{i=1}^{N_{\mathbf{E}}}$. In the spatial discretization of  (\ref{eq:kinetic}), a $\mathcal{Q}^K$-DG scheme seeks a solution in the finite element space
\begin{equation}
    \mathrm{U}_{h}^{K}(\mathcal{T}_{h}) := \left\{ u(\bx): u(\bx)|_{\mathcal{T}_i}\in \mathcal{Q}^{K}(\mathcal{T}_{i}),\;\;1\leq i\leq N_{\mathbf{E}}\right\},
    \label{eq:dg_space}
\end{equation}
where $\mathcal{Q}^{K}(\mathcal{T}_{i})$ is the bi-variate polynomials space defined on the element $\mathcal{T}_i$ whose degree in each direction is at most $K$. Denote the set of cell edges as $\partial\mathcal{T}_{h}$ and the set of edges on the inflow boundary for $\bv_{j}$ as

\begin{equation*}
    \partial\mathcal{T}_{h,j}^{bc} = \left\{ \mathcal{E}:\mathcal{E}\in\partial\mathcal{T}_{h},\mathcal{E}\subset\partial\Omega_{\bx},\bv_{j}\cdot\mathbf{n}(\bx)<0,\;\forall \bx\in\mathcal{E}\right\},
\end{equation*}
where $\mathbf{n}(\bx)$ is the outward normal direction of $\Omega_x$ at $\bx$.   We seek $f_{h}(\bx,\bv_{j},t^n)=f_h^n(\bx,\bv_{j})\in\mathrm{U}_{h}^{K}(\mathcal{T}_{h}), j =1,2,...,N_{\bv}$ such that $\forall\phi_{h}(\bx)\in\mathrm{U}_{h}^{K}(\mathcal{T}_{h})$,
\begin{align}
&\frac{1}{\Dt}\sum_{i=1}^{N_{\mathbf{E}}}\int_{\mathcal{T}_i}\phi_h(\bx) (f_h^n(\bx,\bv_j)-f_h^{n-1}(\bx,\bv_j))d\bx -\sum_{i=1}^{N_{\mathbf{E}}}\int_{\mathcal{T}_i} \Big(\bv_j\cdot\nabla\phi_h(\bx)\Big) f_h^n(\bx,\bv_j) d\bx\notag\\
&+\sum_{\mathcal{E}\in\partial\mathcal{T}_h\setminus\partial\mathcal{T}_{h,j}^{\textrm{bc}}} \int_{\mathcal{E}} \widehat{\BH}(\bv_j,f^n_h, \bn(\bx))\phi_h(\bx) d\bx \notag 
+\sum_{i=1}^{N_{\mathbf{E}}}\int_{\mathcal{T}_i}\sigma_t(\bx) f_h^n(\bx,\bv_j)\phi_h(\bx) d\bx\notag\\
&=\sum_{i=1}^{N_{\mathbf{E}}}\int_{\mathcal{T}_i}(\sigma_s(\bx) \rho_h^n(\bx)+G(\bx))\phi_h(\bx) d\bx
-\sum_{\mathcal{E}\in\partial\mathcal{T}_{h,j}^{\textrm{bc}}}\int_{\mathcal{E}} g(\bx)\phi_h(\bx) \bv_j\cdot \bn(\bx) d\bx,
 \label{eq:DG}
\end{align}
where $\rho_h^n=\sum_{j=1}^{N_{\bv}}\omega_jf^n_h(\bx,\bv_j)$.
The upwind flux $\widehat{\BH}(\bv_j,\eta_h, \bn(\bx))$ along the interface $\mathcal{E}$ between an upwind element $\mathcal{T}^{-}$ and its downwind neighbor $\mathcal{T}^{+}$, is defined as
\begin{align*}
\widehat{\BH}(\bv_j, \eta_h,\bn(\bx))\Big|_{\mathcal{E}} = \frac{\bv_j\cdot\bn(\bx)}{2}\left(\eta_h^-(\bx,\bv_j)+\eta_h^+(\bx,\bv_j)\right)+\frac{|\bv_j\cdot \bn(\bx)|}{2}\Big(\eta_h^-(\bx,\bv_j)-\eta_h^+(\bx,\bv_j)\Big)\Big|_{x\in\mathcal{E}},
\label{eq:upwind}
\end{align*}
Here, $\eta_h^{\pm}$ is the restriction of $\eta_h$ to $\mathcal{E}$ from the element $\mathcal{T}^\pm$.

\subsection{Matrix-vector formulation\label{sec:mat-vec}}
Let $\{\phi_i(\bx)\}_{i=1}^{N_{\bx}}$ be an orthonormal basis for $\mathrm{U}_{h}^{K}(\mathcal{T}_{h})$. Then $f^n_h(\bx,\bv_j)$ and $\rho_h^n(\bx)$ can be expressed by its  degrees of freedom $\bff_j^n=(f^n_{j,1},\cdots,f^n_{j,N_{\bx}})^{\rm T}\in\mathbb{R}^{N_{\bx}}$ as
\begin{equation*} 
f_h(\bx,\bv_j,t^n)=\sum_{k=1}^{N_{\bx}}f^n_{j,k}\phi_k(\bx), {\;\rm and\;}
  \rho_h(\bx,t^n) =\sum_{k=1}^{N_{\bx}}\rho^n_{k}\phi_k(\bx),{\;\rm with\;}\rho^n_{k}=\sum_{j=1}^{N_{\bv}}\omega_jf^n_{j,k}.
\end{equation*}
Let $\brho^n=(\rho^n_{1},\cdots,\rho^n_{N_{\bx}})^{\rm T}\in\mathbb{R}^{N_{\bx}}$, and then the DG scheme (\ref{eq:DG}) can  be written in its matrix form:
\begin{subequations}
\label{eq:kinetic_dg_discrete}
\begin{align}
\big(\frac{1}{\Delta t}\BM+\BD_j+\BSigma_t\big)\bff_j^{n}&=\BSigma_s\brho^{n}+\frac{1}{\Dt}\BM\bff_j^{n-1}+\bG+\bg_j^{(\textrm{bc})}\triangleq\BSigma_s\brho^{n}+\frac{1}{\Dt}\BM\bff_j^{n-1}+ \widetilde{\bG}_j,
\label{eq:dg_discrete1}\\
\brho^{n}&=\sum_{j=1}^{N_{\bv}}\omega_j\bff_j^{n},\;\;j=1,\dots,N_{\bv}.
\label{eq:dg_discrete2}
\end{align}
\label{dg_scheme}
\end{subequations}

Here, the discrete mass matrix $\BM\in\mathbb{R}^{N_{\bx}\times N_{\bx}}$, the discrete advection operator $\BD_j\in\mathbb{R}^{N_{\bx}\times N_{\bx}}$, the discrete scattering cross section $\BSigma_s\in\mathbb{R}^{N_{\bx}\times N_{\bx}}$, total cross section $\BSigma_t\in\mathbb{R}^{N_{\bx}\times N_{\bx}}$,  the discrete source $\bG\in\mathbb{R}^{N_{\bx}}$, and boundary flux $\bg_j^{(\textrm{bc})}\in\mathbb{R}^{N_{\bx}}$ are defined as: 
\begin{subequations}
\begin{align*}
&(\BM)_{kl} = \sum_{i=1}^{N_{\mathbf{E}}}\int_{\mathcal{T}_i} \phi_k(\bx) \phi_l(\bx) d\bx,\\
    &(\BD_{j})_{kl} = -\sum_{i=1}^{N_{\mathbf{E}}}\int_{\mathcal{T}_i} (\bv_j\cdot\nabla\phi_k(\bx)) \phi_l(\bx) d\bx+ \sum_{\mathcal{E}\in\partial\mathcal{T}_h\setminus\partial\mathcal{T}_{h,j}^{\textrm{bc}}}\int_{\mathcal{E}} \widehat{\BH}\left(\bv_j, \phi_l,\bn(\bx)\right)\phi_k(\bx) d\bx,
\\
&(\BSigma_s)_{kl} = \sum_{i=1}^{N_{\mathbf{E}}}\int_{\mathcal{T}_i} \sigma_s(\bx)\phi_k(\bx) \phi_l(\bx) d\bx,\quad (\BSigma_t)_{kl} = \sum_{i=1}^{N_{\mathbf{E}}}\int_{\mathcal{T}_i} \sigma_t(\bx)\phi_k(\bx) \phi_l(\bx) d\bx,
\label{eq:element_sigma}
\\
&(\bG)_k= \sum_{i=1}^{N_{\mathbf{E}}}\int_{\mathcal{T}_i} G(\bx)\phi_k(\bx) d\bx,\quad
(\bg_j^{(\textrm{bc})})_k = -\sum_{\mathcal{E}\in\partial\mathcal{T}_{h,j}^{\textrm{bc}}}\int_{\mathcal{E}} g(\bx)\phi_k(\bx) \bv_j\cdot \bn(\bx) d\bx.
\end{align*}
\end{subequations}

\subsection{Source iteration with diffusion synthetic acceleration (SI-DSA)\label{sec:sisa}}
In this subsection, we review the basic ideas of Source Iteration with Diffusion Synthetic Acceleration (SI-DSA) \cite{adams2002fast}.  For each time step $t^n$, given an initial guess for the macroscopic density at $t^n=n\Delta t$, namely $\brho^{n,(0)}$, SI-DSA seeks the numerical solution $\bff_j^{n}$ and $\brho^{n}$ in an iteratively. The $l$-th iteration of SI-SA consists of an SI step and an SA step.

\textbf{Source Iteration (SI):}
Freeze the density term in equation \eqref{eq:dg_discrete1} as $\brho^{n,(l-1)}$ and then obtain $\bff_j^{n,(l)}$ by solving
\begin{equation}
    \big(\frac{1}{\Delta t}\BM+\BD_j+\BSigma_t\big)\bff_j^{n,(l)}=\BSigma_s\brho^{n,(l-1)}+\frac{1}{\Dt}\BM\bff_j^{n-1}+ \widetilde{\bG}_j,\;\;j=1,\dots,N_{\bv}.
    \label{eq:discrete_SI_f}
\end{equation}
Subsequently, the density is updated correspondingly as $\brho^{n,(l-\half)}=\sum_{j=1}^{N_{\bv}}\omega_j\bff_j^{n,(l)}$. 

In practice, equation \eqref{eq:discrete_SI_f} can be solved highly efficiently through matrix-free transport sweeps \cite{adams2002fast}. In transport sweeps, spatial elements are reordered along the upwind direction of each angle $\bv_j$ to permute the matrix $\frac{1}{\Delta t}\BM+\BD_j+\BSigma_t$ into a block lower triangular one, which can be solved by a single pass of matrix-free block Gauss-Seidel iteration.

\textbf{Source Acceleration (SA):}
If one simply set $\brho^{n,(l)}=\brho^{n,(l-\half)}=\sum_{j=1}^{N_{\bv}}\omega_j\bff_j^{n,(l)}$, SI may suffer from arbitrarily slow convergence for optically thick problems    \cite{adams2002fast}. SA can be understood as a preconditioner to resolve this issue by introducing a correction to the density:
$\brho^{n,(l)}=\brho^{n,(l-\half)}+\delta\brho^{n,(l)}.\label{eq:ideal_density}$
The ideal density correction that ensures convergence in the next iteration is:
\begin{equation}
    \delta\brho^{n,(l)} = \brho^n-\brho^{n,(l-\half)}=\sum_{j=1}^{N_{\bv}}\omega_j\delta\bff_j^{n,(l)},\quad\text{where } \delta \bff_{j}^{n,(l)} = \bff_j^{n} - \bff_j^{n,(l)}.
\label{eq:ideal_density}
\end{equation}
Subtracting Eq. (\ref{eq:discrete_SI_f}) from Eq. (\ref{eq:dg_discrete1}) yields a discrete kinetic correction equation:
\begin{equation}
     \big(\frac{1}{\Delta t}\BM+\BD_j+\BSigma_t\big) \delta\bff_j^{n,(l)} =\BSigma_s\delta\brho^{n,(l)} + \BSigma_s(\brho^{n,(l-\half)}-\brho^{n,(l-1)}),\;\delta\brho^{n,(l)}= \sum_{j=1}^{N_{\bv}}\omega_j\delta\bff_j^{n,(l)},\; 1\leq j\leq N_{\bv},
\label{eq:discrete_correction_f}
\end{equation}
This correction equation  can be seen as a discretization to a steady state RTE with a modified absorption cross section $\tilde{\sigma}_a=\frac{1}{\Delta t}+\sigma_a$, isotropic source $\sigma_s(\rho^{n,(l-\half)}-\rho^{n,(l-1)})$ and zero inflow boundary conditions.

However, the computational cost of solving the ideal correction equation \eqref{eq:discrete_correction_f} is almost the same as solving the original RTE \eqref{eq:rte}. Hence, in DSA \cite{kopp1963synthetic,alcouffe1977diffusion,adams1992diffusion,wareing1993new}, we solve its diffusion limit instead: 
\begin{equation}
\mathbf{C}\delta\brho^{n,(l)}= \BSigma_s(\brho^{n,(l-\half)}-\brho^{n,(l-1)}).
   \label{eq:rho_residual}
\end{equation}
Alternative low-order approximations can also be applied, e.g.  low-order $S_N$ angular dsicretization in 
${\rm S}_2$-SA \cite{lorence1989s} and transport SA \cite{ramone1997transport}, and a second order moment closure in quasi-diffusion  method \cite{gol1964quasi,anistratov1993nonlinear,olivier2023family}.

\textbf{Memory efficient formulation of SI-SA:} In practice, a more memory-efficient formulation that iteratively update only $\rho$ can be applied and derived as follows. Multiplying $\big(\frac{1}{\Delta t}\BM+\BD_j+\BSigma_t\big)^{-1}$
to both sides of equation \eqref{eq:dg_discrete1} followed by a numerical integration in the angular space, we obtain 
\begin{align}
\brho^n=\BT\BSigma_s\brho^n+\wbb^{n-1}
\Longleftrightarrow
    \WBA\brho^n\triangleq(\BI-\BT\BSigma_s)\brho^n=\wbb^{n-1},
    \label{eq:implicit_linear_system_rho}
\end{align}
where $\BT = \sum_{j=1}^{N_{\bv}} \omega_j\big(\frac{1}{\Dt}\BM+\BD_j+\BSigma_t\big)^{-1}$ and  $\widetilde{\bb}^{n-1}=\sum_{j=1}^{N_{\bv}}\omega_j\big(\frac{1}{\Dt}\BM+\BD_j+\BSigma_t\big)^{-1}(\frac{1}{\Dt}\BM\bff_j^{n-1}+\widetilde{\bG}_j)$. Computationally, the matrix-vector product determined by $\BT$ is implemented through matrix-free transport sweeps and numerical integration, requiring $O(N_{\bx}N_{\bv})$ cost.
 Similarly, the SI equation \eqref{eq:discrete_SI_f}  can be reformulated into 
\begin{align}
    &\brho^{n,(l-\half)} = \BT\BSigma_s\brho^{n,(l-1)}+\widetilde{\bb}^{n-1}.
    \label{eq:SI_rho}
\end{align}
The corresponding ideal correction equation can be rewritten as:
\begin{align}
    &(\BI-\BT\BSigma_s)\delta\brho^{n,(l)} =  \BT\BSigma_s(\brho^{n,(l-\half)}-\brho^{n,(l-1)})\notag\\
    \Longleftrightarrow
    &\WBC\delta\brho^{n,(l)}\triangleq (\BT\BSigma_s)^{-1}(\BI-\BT\BSigma_s)\delta\brho^{n,(l)}=\brho^{n,(l-\half)}-\brho^{n,(l-1)}.
\label{eq_ideal_correct}
\end{align}

Note that $\wbb^{n}$ relies on $\bff_j^n$. Hence, after the convergence of SI-DSA to compute $\brho^{n}$, additional transport sweeps are applied to obtain $\bff^n_j$ by solving 
\begin{equation}
(\frac{1}{\Delta t}\BM+\BD_j+\BSigma_t)\bff_j^n=\BSigma_s\brho^n+\wbG_j,\;j=1,\dots,N_{\bv}. 
\label{eq:sweep_to_obtain_f}
\end{equation}
The memory-efficient formulation of SI-DSA is summarized in Alg. \ref{alg:sisa}.

In summary, in each time step, we need one transport sweeps to compute the right hand side of the memory-efficient formulation $\wbb^{n-1}$, multiple transport sweeps to solve $(\BI-\BT\BSigma_s)\brho^n=\wbb^{n-1}$ via SI-SA, and one additional transport sweep to obtain $\bff_j^n$ from $\brho^n$.

\subsection{Motivation of our work \label{sec:motivation}}

Despite the success of classical SI-SA such as SI-DSA, several limitations persist. Their efficiency relies on the validity of underlying empirical assumptions. For example, DSA may become less effective when the kinetic equation is not well approximated by its diffusion limit \cite{ren2019fast}. In addition, when the solution varies smoothly in time, the solution manifold often exhibits low-rank structures across time. However, classical SI-SA fail to exploit such structures.

For steady-state parametric problems, reduced order models (ROMs) have been exploited to construct enhanced SA preconditioners. Instead of relying on empirical low-order approximations such as the diffusion limit, ROMs are constructed directly from the original kinetic formulation by extracting low-rank structures from solution data across parameters \cite{peng2024reduced,peng2025flexible,tang2025synthetic}. Specifically, an offline-online decomposition framework is employed: in the offline stage, ROMs are constructed by  learning the low-rank structures of solution manifold from pre-collected training solution data; in the online stage, these ROMs are used to provide improved initial guesses and enhanced SA corrections.

However, in non-parametric long-time simulations, the offline-online decomposition is no longer applicable, since solution data are generated on-the-fly during time marching rather than available in advance. In addition, the approaches in \cite{peng2024reduced,peng2025flexible,tang2025synthetic} construct ROMs using data for the microscopic particle distribution, which leads to high memory cost. Although a projection-based ROM for the steady-state counterpart of the memory-efficient reformulation \eqref{eq:implicit_linear_system_rho} has been proposed in \cite{behne2022minimally}, projecting $\BI - \BT \BSigma_s$ requires $\mathcal{O}(N_{\bx}N_{\bv}r^2)$ computational cost with $r$ denoting the rank.

In this paper, we aim to address these challenges by designing an efficient and memory-lean on-the-fly ROM-based acceleration for SI-DSA in implicit time marching.

\begin{algorithm}[H]
\caption{Memory-efficient Source Iteration and Synthetic Acceleration. \label{alg:sisa} }
\begin{algorithmic}[1]
\STATE{Given an initial guess $\brho^{n,(0)}$, tolerance $\epsilon_{\mathrm{SISA}}$, the maximum number of iterations $N_{\mathrm{iter}}$, and low-rank operator $\WBC$ for correction equation.}
\FOR{ $l = 1:N_{\mathrm{iter}}$}
    \STATE{\textbf{Source Iteration:}}
    \STATE{Solve $\brho^{n,(l-\half)} = \BT\BSigma_s\brho^{n,(l-1)}+\widetilde{\bb}^{n-1}$ with transport sweep.} 
    \IF{ $\|\brho^{n,(l-\half)} - \brho^{n,(l-1))}\|_{\infty} < \epsilon_{\mathrm{SISA}}$}
        \STATE{return solutions: $\brho^{n,(l-\half)}$ and $\bff^{n,(l)}$}.
    \ELSE{
    \STATE{\textbf{Synthetic acceleration:}}
    \STATE{Solve the correction equation $\WBC \delta\brho^{n,(l)}=\brho^{n,(l-\half)}-\brho^{n,(l-1)}$ to obtain the correction $\delta \brho^{n,(l)}$}.
    \STATE{Update the density flux: $\brho^{n,(l)} = \brho^{n,(l-\half)} + \delta \brho^{n,(l)}$}.}
    \ENDIF
\ENDFOR
\end{algorithmic}
\end{algorithm}
\section{On-the-fly reduced order model (ROM)-based acceleration for time-dependent RTE}
\label{sec:ROM_construct}
In this paper, we aim to  leverage ROMs built on-the-fly during implicit time marching of RTE to accelerate the underlying SI-DSA solver. 

To extract low-rank structures across time from solution  data generated on-the-fly during time marching, we leverage dynamical mode decomposition (DMD) and incremental singular value decomposition (SVD) to construct and adaptively update ROMs to enhance both initial guesses and preconditioners.  
Moreover, the ROM is directly constructed for the memory-efficient reformulation \eqref{eq:implicit_linear_system_rho} in a purely data-driven manner. In other words, high memory costs in \cite{peng2022reduced,peng2024micro,tang2025synthetic} and high computational cost due to projections in \cite{behne2022minimally} are avoided at the same time. 

In this section, we first present  DMD algorithm to construct ROMs in Sec. \ref{sec:dmd}, then outline our ROM-based accelerations for SI-DSA and adaptive update of ROMs in Sec. \ref{sec:rom-sisa}. In in Sec. \ref{sec:comparison}, to contextualize our method, we compare it with the DMD-based preconditioner in \cite{mcclarren2022data} and point out their complementary power.

\subsection{Dynamic mode decomposition with incremental SVD\label{sec:dmd}}
Here, we present dynamical mode decomposition (DMD) \cite{schmid2011applications,tu2013dynamic,schmid2011application} used to construct ROMs in a purely data-driven manner. Particularly, we follow the variable DMD in \cite{smith2023variable}. At time step $n$, we solve the linear system $\WBA\brho^n=\wbb^{n-1}$ (defined in \eqref{eq:implicit_linear_system_rho}) in implicit time marching.
We construct two data/snapshot matrices collecting the solution and the corresponding right handside at time $t_n$:
\begin{equation}
    \mathbf{R}_m = \begin{bmatrix}
        \brho^{1}, &\brho^{2}, &\cdots,& \brho^{m}
    \end{bmatrix}\in\mathbb{R}^{N_{\bx}\times m} \quad\text{and}\quad
    \BB_m = \begin{bmatrix}
        \wbb^{0}, &\wbb^{1}, &\cdots,& \wbb^{m-1}
    \end{bmatrix}\in\mathbb{R}^{N_{\bx}\times m}.
    \label{eq:snapshots}
\end{equation}
Then, $\WBA\BR_m=\BB_m$.
Compute a truncated SVD of the matrix $\BR_m\approx\BU_{m}\BS_{m}\BV_{m}^{\rT}$, 
where $\BU_{m}\in\mathbb{R}^{N_{\bx}\times r_m},\BS_{m}=\textrm{diag}(s_1,\dots,s_{r_m})\in\mathbb{R}^{r_m\times r_m},\BV_{m}\in\mathbb{R}^{m\times r_m}$, and $r_m\leq\min\{m,N_{\bx}\}$ is the rank.
Then, we can derive a reduced system, namely $\WBA_{r_m}$, that approximates the Galerkin projection of $\WBA$ onto the linear space spanned by $\BU_{m}$:
\begin{equation}
    \WBA\BU_{m}\BS_{m}\BV_{m}^\rT\approx \BB_m\Rightarrow
    (\BU_{m})^\rT\WBA\BU_{m}\approx \BU_{m}^\rT\BB_{m}\BV_{m}\BS_{m}^{-1}\triangleq\WBA_{r_m}
    \label{eq:dmd-operator}
\end{equation}
Then, $\brho^n$ can be approximated by
\begin{equation}
    \brho^{n}\approx \BU_{m}\bc^n,\;\text{where }\bc^n\;\text{solves}\;\WBA_{r_m}\bc^n=\BU_{m}^\rT\wbb^{n-1}.
\end{equation}

\textbf{Incremental SVD for streaming data.} In the crucial step to compute the SVD of $\BR_m$, we apply the incremental SVD through low-rank modification in \cite{brand2006fast}. For brevity, details of this algorithm is summarized in \ref{appx:i-svd}. 
Given $\BR_{m-1}\approx \BU_{m-1}\BS_{m-1}\BV_{m-1}^\rT$, this algorithm builds the SVD for $\BR_{m}=\begin{bmatrix}\BR_{m-1}&\brho^m\end{bmatrix}$ very efficiently through projections and SVD of a low-dimensional core matrix whose size is $r_{m-1}+1$. By applying this algorithm in each time step, the SVD of $\BR_m$ can be built incrementally as time marching proceed.  The main benefit of this strategy include:
\begin{enumerate}
    \item The solution data is generating in a streaming manner. When leveraging DMD to construct enhanced initial guesses and preconditioners, it is difficult to determine a priori how many data snapshots are required to obtain a sufficiently accurate ROM. Incremental update allows us to dynamically determine the size of required data on-the-fly through error indicators.
    \item Given the rank-$r$ SVD of $\BR_{m-1}\in\mathbb{R}^{N_{\bx}\times m}$, the incremental SVD update  in \cite{brand2006fast} computes the SVD of $O((N_{\bx}+m+r_m)r^2_m)$ cost. This strategy is especially efficient when ROM needs to be adaptively updated.
\end{enumerate}

\subsection{On-the-fly ROM-based acceleration of SI-DSA for implicit time marching \label{sec:rom-sisa}}
\begin{figure}[t]
\centering
\begin{tikzpicture}[
    >=Latex,
    node distance=25mm,
    box/.style={
        draw,
        rounded corners=2mm,
        align=center,
        minimum width=3.8cm,
        minimum height=1.2cm
    },
    arrow/.style={->, thick}
]

\node[box] (p1) {Phase 1:\\ Construct ROM for IG};

\node[box, right=of p1] (p2) {Phase 2:\\ Construct ROM for PC};

\node[box, right=of p2] (p3) {Phase 3:\\ Apply ROM-IG and ROM-PC\\ Adaptively update ROMs};

\draw[arrow] (p1) -- node[above] {\small{Apply ROM-IG}} (p2);
\draw[arrow] (p2) -- (p3);

\end{tikzpicture}
\caption{Flowchart for On-the-fly ROM-based acceleration workflow. IG: initial guesses. PC: preconditioner. ROM-IG/PC: ROM-enhanced initial guesses/preconditioner.}
\label{fig:flow_chart}
\end{figure}
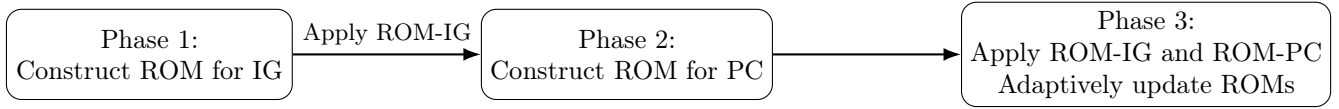

Our on-the-fly ROM-based acceleration strategy includes three phases illustrated in the flowchart Fig. \ref{fig:flow_chart}. Specifically, in phase $1$, we construct  a ROM for the linear system in implicit time marching, i.e. $\WBA\brho^n=\wbb^{n-1}$, to provide improved initial guesses. In phase $2$, we apply ROM-enhanced initial guesses and construct a ROM for the ideal correction equation to enhance preconditioners. In phase $3$, we use ROMs to enhance both initial guesses and preconditioners to accelerate SI-DSA, and adaptively update these ROMs to maintain their accuracy. Below, we elaborate details of each phase.

\textbf{Phase 1: Construct a ROM to provide improved initial guesses.} Following Sec. \ref{sec:dmd} and \ref{appx:i-svd}, we employ  incremental SVD without truncation to iteratively construct a ROM for the discrete system $\WBA\brho^n=\wbb^{n-1}$ via DMD  during time marching until the first time step $N_0$ satisfying 
\begin{equation}
s_{N_0}/{\sum_{k=1}^{N_0}s_k}\leq\epsilon_{\rm IG}.
\label{eq:ig_stop}
\end{equation}
Here, $\epsilon_{\rm IG}$ is a threshold and $s_k$ is the $k$-th singular value for $\BR_{N_0}$ defined in \eqref{eq:snapshots}. Let $\BR_{N_0}=\BU_{\textrm{IG}}\BS_{\textrm{IG}}\BV_{\textrm{IG}}^\rT$, then the corresponding reduced system is $(\BU_{\textrm{IG}})^{\rT}\WBA\BU_{\textrm{IG}}\approx\WBA_{r}
=\BU_{\textrm{IG}}^{\rT}\BB_{N_0}\BV_{\textrm{IG}}\BS_{\textrm{IG}}^{-1}$. For time step $n>N_0$, we use this ROM to provide an improved initial guess as
\begin{equation}
 \brho_{\textrm{DMD}}^{n,(0)}=\BU_{\textrm{IG}}\bc^n,\;\text{where } \WBA_{r}\bc^n=\BU_{\textrm{IG}}^\rT\wbb^{n-1}.
 \label{eq:dmd-ig}
\end{equation}

\textbf{Phase 2: Construct a ROM to enhance the SA preconditioner and adaptively update the ROM for initial-guess.}
Once the ROM for the initial guess is available, we apply the ROM-enhanced initial guess in \eqref{eq:dmd-ig} and proceed to construct a ROM for the ideal correction equation of the first iteration (defined in \eqref{eq_ideal_correct})
\[
\WBC\delta\brho^{n,(1)}=\brho^{n,(\half)}-\brho^{n,(0)}, \qquad n\ge N_0+1.
\]
This ROM will later be used to build a ROM-enhanced preconditioner.
Specifically, we collect data matrices
\begin{subequations}
\begin{align}
    \delta\BR_{m}^{(1)}&=[\delta\brho^{N_0+1,(1)},\cdots,\delta\brho^{N_0+m,(1)}]=[\brho^{N_0+1}-\brho^{N_0+1,(\half)},\cdots,\brho^{N_0+m}-\brho^{N_0+m,(\half)}],\\
    \text{and}\quad
    \delta\BB_m^{(1)}&= [\brho^{N_0+1,(\half)}-\brho^{N_0+1,(0)}_\textrm{DMD},\cdots,\brho^{N_0+m,(\half)}-\brho^{N_0+m,(0)}_\textrm{DMD}]
\end{align}
\label{eq:correction_snapshot}
\end{subequations}
in a streaming manner and  update the DMD-ROM via incremental SVD  following Sec. \ref{sec:dmd}. Assume $\delta\BR_m^{(1)}=\delta\BU_m^{(1)}\delta\BS_m^{(1)}\left(\delta\BV_m^{(1)}\right)^\rT$, where $\delta\BS_m^{(1)}=\textrm{diag}(\delta s_{1}^{(1)},\dots,\delta s_{m}^{(1)})$. We keep incrementally updating the ROM for the correction equation without truncation  via DMD until $m=N_1$ satisfying
\begin{equation}
  \delta s^{(1)}_{N_1}/\sum_{i=1}^{N_1} \delta s^{(1)}_{i}\leq\epsilon_{\rm PC}, 
  \label{eq:pc_stop}
\end{equation} where $\epsilon_{\rm PC}$ is a threshold. Then, we compute the  SVD 
$\delta\BR_{\textrm{PC}}^{(1)}\approx\delta\BU_{\textrm{PC}}^{(1)}\delta\BS_{\textrm{PC}}^{(1)}\left(\delta\BV_{\textrm{PC}}^{(1)}\right)^\rT$.
The reduced order system given by DMD is defined as
\begin{equation}
\left(\delta\BU_{\textrm{PC}}^{(1)}\right)^\rT\WBC\delta\BU_{\textrm{PC}}^{(1)} \approx  \WBC^{(1)}_{r}
\triangleq\left(\delta\BU_{\textrm{PC}}^{(1)}\right)^\rT\delta\BB_{N_1}^{(1)}\delta\BV_{\textrm{PC}}^{(1)}(\delta\BS_{\textrm{PC}}^{(1)})^{-1}.
\end{equation}
This ROM for the correction equation can be leveraged to enhance the SA preconditioning step.

Besides building the ROM for the correction equation, we also simultaneously monitor the accuracy of the ROM for the initial guess to adaptively update it. The adaptation strategy will be deferred to the presentation of phase 3.

\textbf{Phase 3: Apply ROM-enhanced initial guesses and preconditioners, and adaptively update ROMs.} 
When both ROMs for the initial guess and preconditioner are available i.e. $n\geq N_0+N_1+1$, we apply a ROM-enhanced initial guess defined in \eqref{eq:dmd-ig} and a ROM-enhanced hybrid SA preconditioner. In the $l$-th iteration, the correction in the SA step are constructed as follows:
\begin{enumerate}
    \item When $l=1$, the density is correted by $\brho^{n,(1)}=\brho^{n,(\half)}+\delta\brho^{n,(1)}_{\textrm{DMD}}$ with $\delta\brho^{n,(1)}_{\textrm{DMD}}=\delta\BU_{\textrm{PC}}^{(1)}\delta\bc^{n,(1)}$, where $\delta\bc^{n,(1)}$ solves
    \begin{equation}
        \WBC^{(1)}_{r}\delta\bc^{n,(1)}=\left(\delta\BU_{\textrm{PC}}^{(1)}\right)^\rT(\brho^{n,(\half)}-\brho^{n,(0)}_{\textrm{DMD}}).
        \label{eq:dmd-sa}
    \end{equation}
    \item When $l\geq 2$, the density correction is obtained by solving $\BC_{\textrm{DSA}}\delta\brho^{n,(l)}=\BSigma_s(\brho^{n,(l-\half)}-\brho^{n,(l-1)})$, where $\BC_{\textrm{DSA}}$ is the diffusion system used in DSA.
\end{enumerate}
We note that solving equation \eqref{eq:dmd-sa} is computationally cheap since its size equals to the rank of the ROM.

Simultaneously, we monitor the accuracy of constructed ROMs and adaptively update them if necessary.  Specifically, given an accuracy threshold $\epsilon_{\textrm{UP}}$, if
$\|\brho^{n,(0)}_{\textrm{DMD}}-\brho^n\|_2 > \epsilon_{\textrm{UP}},
$
we update the ROM for the initial guess by augmenting data matrices with $\brho^n$ and $\wbb^{n}$ and applying incremental SVD with truncation to update the DMD-ROM. 
Moreover, if either the initial-guess ROM is updated or
$\|\delta\brho^{n,(1)}_{\textrm{DMD}}-\delta\brho^{n,(1)}\|_2 > \epsilon_{\textrm{UP}}$
at time $t_n$, we update the ROM for the preconditioner by augmenting the correction snapshots with
$\delta\brho^{n,(1)}=\brho^{n}-\brho^{n,(\frac{1}{2})}$
and $\delta\wbb^{n,(1)}=\brho^{n,(\half)}-\brho^{n,(0)}_\textrm{DMD}$ and applying incremental SVD with truncation. The use of incremental SVD enables adaptive updates of the ROMs, while the truncation in the ROM updates is introduced to control the rank of the ROM.

\textbf{Key features of the acceleration algorithm.} Now, we highlight several crucial details of our algorithm that merit emphasis or clarification.
\begin{enumerate}    
    \item \textbf{Sequential construction of ROMs for initial guesses and preconditioners.}
    The ROM for the initial guess must be constructed prior to the one for the preconditioner. As pointed out in \cite{tang2025synthetic}, the right-hand side of the correction equation at the $l$-th iteration depends implicitly on the initial guess and the preconditioners from all preceding iterations. Therefore, without a sequential construction that respects this sequential dependency, ROM-based preconditioners may suffer a significant efficiency reduction. Although we only construct the ROM for the correction equation at the first iteration, the procedure naturally extends to subsequent iterations by constructing the ROMs for the initial guess and for the first $l-1$ iterations before building the ROM for iteration $l$. 
    
    \item \textbf{Efficient and memory lean ROM construction through DMD.} In ROM-enhanced preconditioenrs for parametric steady-state RTE \cite{tano2021affine,peng2024reduced}, ROMs are directly constructed to the counterpart of the discrete kinetic  equation \eqref{eq:kinetic_dg_discrete} using solution snapshots for the particle distribution $f$, leading to high memory costs.
    
    In addition, as mentioned earlier, 
    $\BT = \sum_{j=1}^{N_{\bv}} \omega_j\Big(\frac{1}{\Dt}\BM+\BD_j+\BSigma_t\Big)^{-1}$
    is implemented through matrix-free transport sweeps for all angular directions. 
    A direct projection of $\BI-\BT\BSigma_s$ onto an $r$-dimensional space would require an  $O(N_{\bx} N_{\bv}r^2)$ computational cost. 
    Such overhead may substantially diminish or even eliminate practical computational savings \cite{behne2022minimally}, particularly when the ROM is updated adaptively. 

    In contrast, our DMD-based ROMs are directly built based on data for the density $\rho$ instead of the distribution function $f$, reducing the memory cost per snapshot from $O(N_{\bx}N_{\bv})$ to $O(N_{\bx})$. Moreover, the reduced system is constructed in a purely data-driven manner with $O(N_{\bx}r^2)$ costs, avoiding the expensive direct projection of the sweeping-based operator $\BI-\BT\BSigma_s$. 
    
    \item \textbf{Motivation of adaptive update.} During time marching, new solution features that are not captured by the historical data  may emerge. For example, particles may encounter material interfaces. Hence, adaptive updates of ROMs are necessary to maintain their accuracy. 
\end{enumerate}

At last, we summarize our on-the-fly DMD enhanced acceleration for implicit time marching of RTE in Alg.~\ref{alg:whole_alg} following {\tt{Matlab}} notations.

\begin{algorithm}[H]
\caption{On-the-fly DMD enhanced acceleration for SI-SA in implicit time marching. \label{alg:whole_alg} }
\begin{algorithmic}[1]
\STATE{Given an initial guess $\brho^{n,(0)}$, tolerance $\epsilon_{\mathrm{SISA}}$, tolerance $\epsilon_{\mathrm{IG}}$, tolerance $\epsilon_{\mathrm{PC}}$, tolerance $\epsilon_{\mathrm{UP}}$, the maximum number of iterations $N_{\mathrm{iter}}$.} 
\STATE{{\bf Initialization: } set $\BR=[\;]$, $\BB=[\;]$, $\delta\BR=[\;]$, $\delta\BB=[\;]$, $\textrm{Flag}_{\textrm{IG}}=\textrm{False}$, $\textrm{Flag}_{\textrm{PC}}=\textrm{False}$.}
\FOR{ $n=1:N_t$ }
    \STATE{\textbf{Compute RHS}: Compute $\wbb^{n-1}$ defined in equation \eqref{eq:implicit_linear_system_rho} via transport 
    sweeps.}
    \IF{$\textrm{Flag}_{\textrm{IG}}$}
        \STATE{Use ROM-enhanced initial guess $\brho^{n,(0)}$ defined in \eqref{eq:dmd-ig}.}
    \ELSE
        \STATE{Set initial guess as $\brho^{n,(0)}=\brho^{n-1}$.}
    \ENDIF
    
    \STATE{\textbf{Solve for $\brho^n$ and $\bff_j^n$.}}
    \IF{$\textrm{Flag}_{\textrm{PC}}$}
        \STATE{Use the ROM-enhanced preconditioner outlined in Sec. \ref{sec:rom-sisa} to compute $\brho^n$.}
    \ELSE
        \STATE{Use SI-DSA to compute $\brho^n$.}
    \ENDIF

    \STATE{Obtain $\bff_j^n$ from $\brho^n$ by applying transport sweeps to solve \eqref{eq:sweep_to_obtain_f}.}

    \IF{Not $\textrm{Flag}_{\textrm{IG}}$}
        \STATE{Update snapshot matrices $\BR=[\BR\;\;\brho^{n}]$, $\BB=[\BB\;\;\wbb^{n-1}]$, and then apply incremental SVD to update the reduced system $\WBA_{r}$ via DMD following Sec. \ref{sec:dmd}.}
        \IF{ \eqref{eq:ig_stop} is satisfied}
            \STATE{$\textrm{Flag}_{\textrm{IG}}=\textrm{True}$.}
        \ENDIF
    \ELSE
        \IF{$||\brho^{n,(0)}_{\rm DMD}-\brho^n||_2>\epsilon_{\rm UP}$.}
        \STATE{Adaptively update the ROM for initial guesses via DMD and incremental SVD.}
        \ENDIF
        \IF{Not $\textrm{Flag}_{\textrm{PC}}$}
        \STATE{Update snapshot matrices $\delta\BR=[\delta\BR\;\;\delta\brho^{n,(1)}]$ and $\delta\BB=[\delta\BB\;\;\brho^{n,(\half)}-\brho^{n,(0)}_\textrm{DMD}]$, and then apply incremental SVD to update the reduced correction system $\WBC_{r}^{(1)}$ via DMD.}

         \IF{\eqref{eq:pc_stop} is satisfied}
            \STATE{$\textrm{Flag}_{\textrm{PC}}=\textrm{True}$.}
        \ENDIF
        \ENDIF
        \IF{$||\brho^{n,(0)}_{\rm DMD}-\brho^n||_2>\epsilon_{\rm UP}$ or $||\delta\brho^{n,(1)}_{\rm DMD}-\delta\brho^{n,(1)}||_2>\epsilon_{\rm UP}$.}
        \STATE{Adaptively update the ROM for preconditioner via DMD and incremental SVD.}
        \ENDIF
    \ENDIF
\ENDFOR
\end{algorithmic}
\end{algorithm}

\subsection{Comparison with DMD-based preconditioner in \cite{mcclarren2022data}}
\label{sec:comparison}
In \cite{mcclarren2022data}, a novel DMD-based preconditioner is proposed for nonlinear thermal radiation requiring non-negative fix for the distribution function. 
This approach and our DMD enhanced preconditioner have different goals leading to different but complimentary methodology. 
\begin{enumerate}
    \item One of our goals is to enhance DSA. In contrast, \cite{mcclarren2019calculating} aims to provide an alternative for DSA when non-negative fixing makes classical DSA preconditioner inapplicable. 
    \item We build ROM by extracting low-rank structures from the solution data of the ideal correction equation across time:
    \begin{equation}
        [\brho^{N_0+1}-\brho^{N_0+1,(\half)},\dots,\brho^{N_0+N_1}-\brho^{N_0+N_1,(\half)}]\quad \text{and}\quad [\delta\wbb^{N_0+1,(1)},\dots,\delta\wbb^{N_0+N_1,(1)}].
    \end{equation}
    In contrast, \cite{mcclarren2022data}
     builds ROM by extracting low-rank structures across source iterations without SA within the same time step. Specifically, they apply DMD to different data matrices:
    \begin{align}    
    &[\brho^{n,(1)}-\brho^{n,(0)}, \brho^{n,(2)}-\brho^{n,(1)}, \cdots, \brho^{n,(m)}-\brho^{n,(m-1)}
    ]\notag\\
    \quad\text{and}\quad 
    &[
        \brho^{n,(2)}-\brho^{n,(1)}, \brho^{n,(3)}-\brho^{n,(2)}, \cdots, \brho^{n,(m+1)}-\brho^{n,(m)}
    ].
    \end{align}
\end{enumerate}
On the other hand, our method and the approach in \cite{mcclarren2019calculating} are complimentary to each other. First, we leverage low-rank structures across time, while they utilize low-rank structures across iterations within the same time step. Moreover, 
as an alternative to DSA, the approach in \cite{mcclarren2019calculating} can be directly integrated into our framework by replacing DSA  with it.

\section{Numerical results}
\label{sec:numerical}
We demonstrate the performance of the proposed on-the-fly ROM-based acceleration for SI-DSA through a series of numerical tests in 1D slab geometry and 2D $X$-$Y$ geometry. 

We implement our code in {\tt{Matlab}}. A fully consistent version of DSA \cite{adams2001discontinuous} is applied, and its detailed formulation can be derived following \ref{appx:i-svd} of \cite{peng2024reduced} by replacing $\BSigma_t$ in it with $\BSigma_t+\frac{1}{\Delta t}\BM$. 
We solve the transport equations in the SI step with matrix-free transport sweeps, the diffusion
  equation in DSA using conjugate gradient method with algebraic multigrid preconditioner implemented in {\tt{iFEM}} package \cite{chen2008ifem}, and direct LU decomposition for the reduced systems due to their small sizes. Unless otherwise specified, we set $\epsilon_{\textrm{SISA}}=10^{-11}$, $\epsilon_{\textrm{IG}}=10^{-9}$, $\epsilon_{\textrm{PC}}=10^{-6}$, and $\epsilon_{\textrm{UP}}=10^{-9}$. In all tests, we use piecewise linear polynomial DG.

We compare our method with SI-DSA, and refer our method as ``DMD enhanced SI-DSA" if DSA is in the ROM-enhanced hybrid preconditioner. We denote solution given by SI-DSA as $\brho^n_{\textrm{SI-DSA}}$ and the solution given by DMD enhanced SI-DSA as $\brho^n_{\textrm{DMD-SI-SA}}$. Particularly, we define 
\begin{subequations}
\begin{align}
    &\textrm{Avg. sweeps per time step}=\frac{\textrm{Total number of transport sweeps for time steps over a period}}{\textrm{Total number of time steps in this preiod}},\\
    &\textrm{Avg. relative computational time}=\frac{\textrm{Computational time for time marching with DMD enhanced SI-DSA}}{\textrm{Computational time for time-marching using SI-DSA}}.
\end{align}
\end{subequations}
Throughout this numerical section, ``IG" refers to initial guess, ``PC" refers to preconditioners, and ``UP" refers to adaptive update. Recall that our on-the-fly ROM acceleration consists three phases: 
Phase I, online construction of a ROM to provide enhanced initial guesses; Phase II, online construction of a ROM-based preconditioner and application of ROM-enhanced initial guesses; and Phase III, application and adaptive updates of the ROM-enhanced initial guesses and ROM-enhanced preconditioners.

\subsection{Two-material problem in 1D slab geometry\label{sec:1d-two-material}}
We consider a two-material problem on $\Omega_x=[0,11]$ and $\Omega_{v}=[-1,1]$ with initial condition $f(x,v,0)=0$ and isotropic inflow boundary conditions:
\begin{subequations}
\begin{align}
&\sigma_s(x)=\left\{\begin{array}{ll}
0, & \text{if}   \;x\in[0,1], \\
100, & \text{if}\;x\in[1,11],\quad
\end{array}\right.\quad
\sigma_a(x)=\left\{\begin{array}{ll}
1, & \text{if}   \;x\in[0,1], \\
0, & \text{if}\;x\in[1,11],\quad
\end{array}\right.\quad \sigma_t(x)=\sigma_a(x)+\sigma_s(x).\\
&f(0,v,t)=5 {\; \rm with \;}\bv>0,\quad\text{and}\quad f(11,v,t)=0 {\;\rm with \;}\bv<0,
\end{align}
\label{eq:discrete_correction}
\end{subequations}
 In this test, \pzc{$6$}-point Gaussian-Legendre quadrature points for [-1,1] are used for angular discretization, and a mesh size with $\Delta x=0.1$ is used to partition the spatial space. 

\begin{figure}[htbp]
    \centering
    \includegraphics[height=0.33\textwidth]{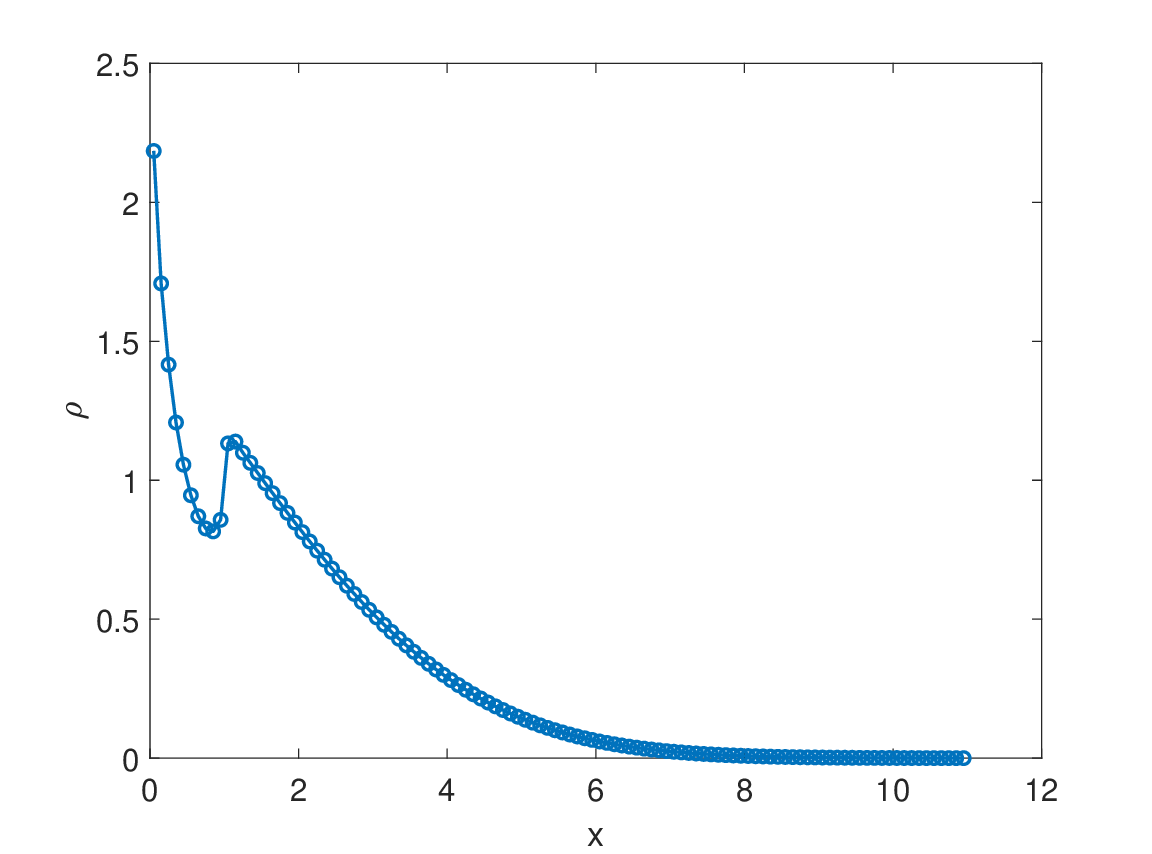} \includegraphics[height=0.33\textwidth]{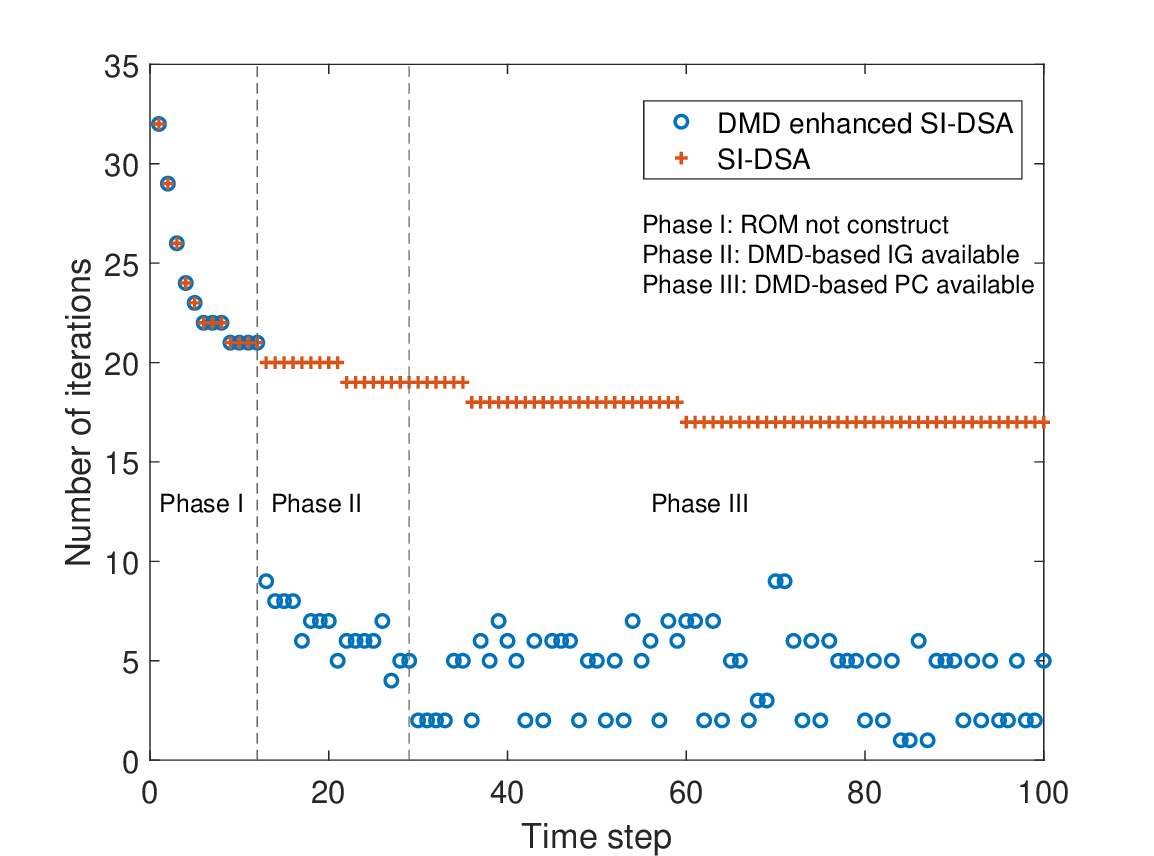} 
    \caption{Results for the two-material problem in Sec. \ref{sec:1d-two-material}. Left: $\rho$  at $T=1000$. Right: number of iterations required to convergence at each time step.\label{fig:two_material_solution}}  
\end{figure}

We solve the problem from $t=0$ to $t=1000$ with time-step size $\Delta t=10$. The solution at $t=1000$ is presented in the left plot of Fig. \ref{fig:two_material_solution}. The difference between the solution vector at $t=1000$ given by SI-DSA and our approach is $||\brho^{100}_{\textrm{SI-DSA}}-\brho^{100}_{\textrm{DMD-SI-SA}}||_2=1.02\times10^{-9}$.  

\textbf{Acceleration.} To evaluate the efficiency of the proposed method, we present the number of iterations for convergence required by SI-DSA and our approach in the right plot of Fig. \ref{fig:two_material_solution}. As the ROM to enhance initial guesses and preconditioners are progressively built, we observe that the iteration count decreases progressively from Phase 1 to Phase 3. In this test, our ROM-based acceleration  effectively reduce the number of iterations for convergence to approximately $2–5$, while SI-DSA requires at least $17$ iterations for convergence and even more than $20$ in early stages.

\subsection{Isotropic scattering\label{sec:iso}}
\begin{figure}[htbp]
    \centering
    \begin{subfigure}[b]{0.33\textwidth}
        \centering
        \hspace{-20mm}
        \begin{minipage}{4cm}
        \includegraphics[scale=0.44]{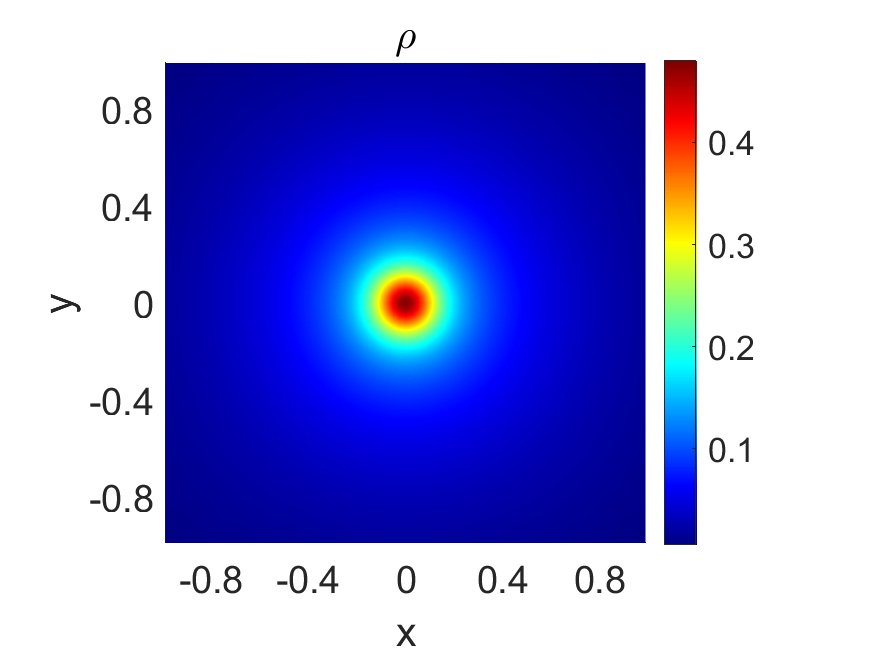}
        \end{minipage}
        \label{fig:scattering_rho}
    \end{subfigure}
     \hspace{-0.01\textwidth}
     \begin{subfigure}[b]{0.33\textwidth}
        \centering
        \hspace{-20mm}
        \begin{minipage}{4cm}
        \includegraphics[scale=0.44]{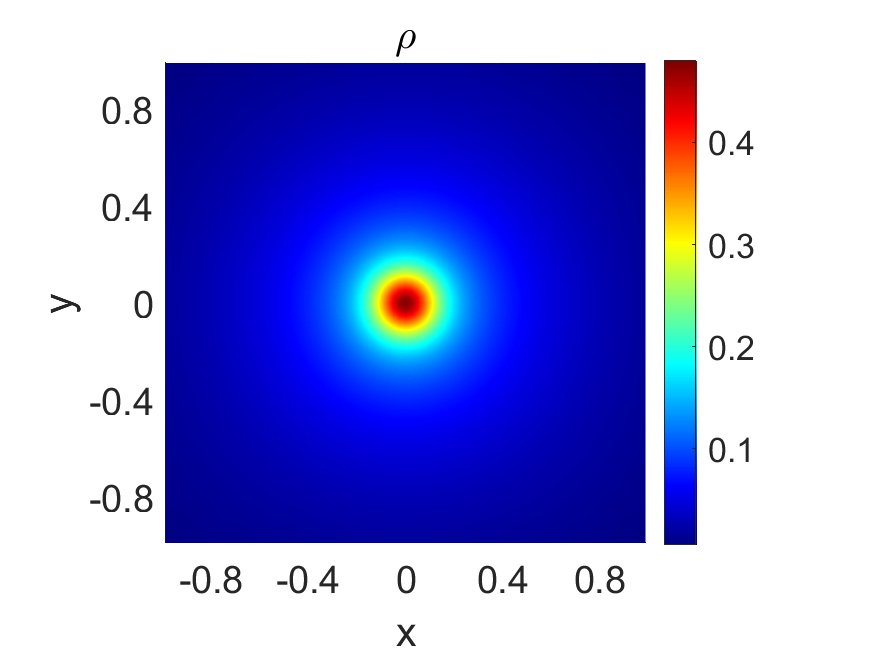}
        \end{minipage}
        \label{fig:scattering_rho}
    \end{subfigure}
     \hspace{-0.01\textwidth}
    \begin{subfigure}[b]{0.33\textwidth}
        \centering
        \hspace{-10mm}
        \begin{minipage}{4cm}
        \includegraphics[scale=0.39]{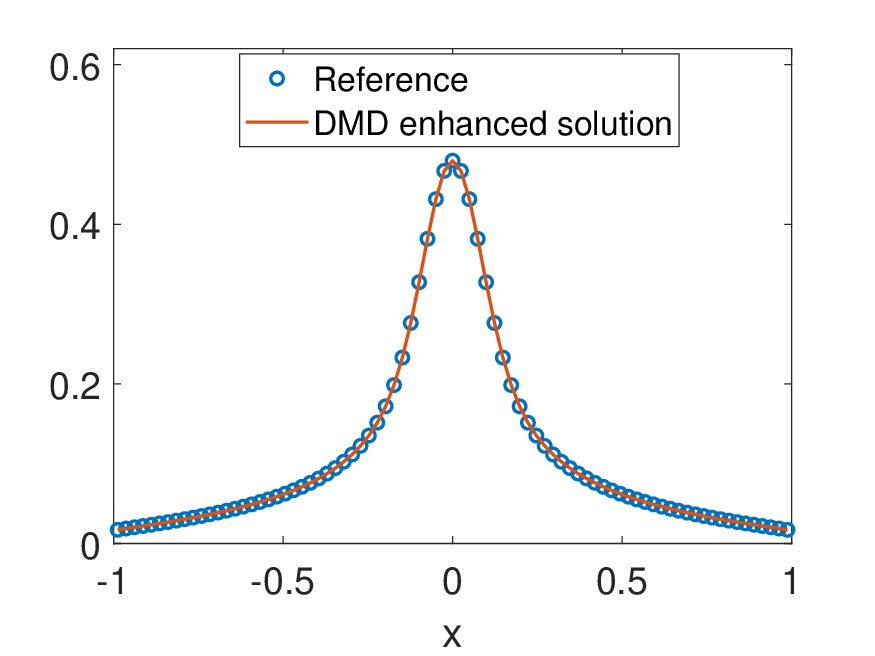}
        \end{minipage}
        \label{fig:scattering_rho_y0}
    \end{subfigure}
    \caption{Density profile at t=2.5 with condition CFL=1 and $\sigma_s=1$. Left: solution by DMD enhanced SI-DSA, where we set $\epsilon_{\textrm{IG}}=10^{-9}$, $\epsilon_{\textrm{PC}}=10^{-6}$, and $\epsilon_{\textrm{UP}}=10^{-9}$. Middle: reference solution by SI-DSA. Right: the comparison of two solutions with $y=0$.}    \label{fig:Gaussian_source_variable}
\end{figure}

We consider a problem on the computational domain $[-1,1]^2$ with isotropic scattering and no absorption ($\sigma_a=0$), zero inflow boundary conditions, zero initial conditions and a Gaussian source,
\begin{equation}
\bG(\bx, \by) = \frac{10}{\pi}{\rm exp}(-100(x^2+y^2)).
\end{equation}
We consider $\sigma_s=0.1,1,10,100$.
We partition the computational domain using a uniform $81\times 81$ rectangular mesh, and use CL($40,6$) quadrature in the angular space. We solve the problem from $t=0$ to $t=2.5$. In this example, the tolerance in the stopping criteria of SISA is set to $\epsilon_{\textrm{SISA}}=10^{-12}$. We present the reference density profile for the density profile at $t=2.5$ for $\sigma_s=1$ in Fig. \ref{fig:Gaussian_source_variable} and compare it to the DMD enhanced SI-DSA solution with $\Delta t=\Delta x$ and $\epsilon_{\textrm{IG}}=\epsilon_{\textrm{UP}}=10^{-9}$ and $\epsilon_{\textrm{PC}}=10^{-6}$. Their difference in the $l_2$ norm is $\|\rho_{\rm DMD} - \rho\|_2 = 9.51 \times 10^{-16}$.

\begin{figure}[h]
    \centering
    \begin{subfigure}[b]{0.45\textwidth}
        \centering
        \begin{minipage}{6cm}
        \hspace{-10mm}
        \includegraphics[scale=0.55]{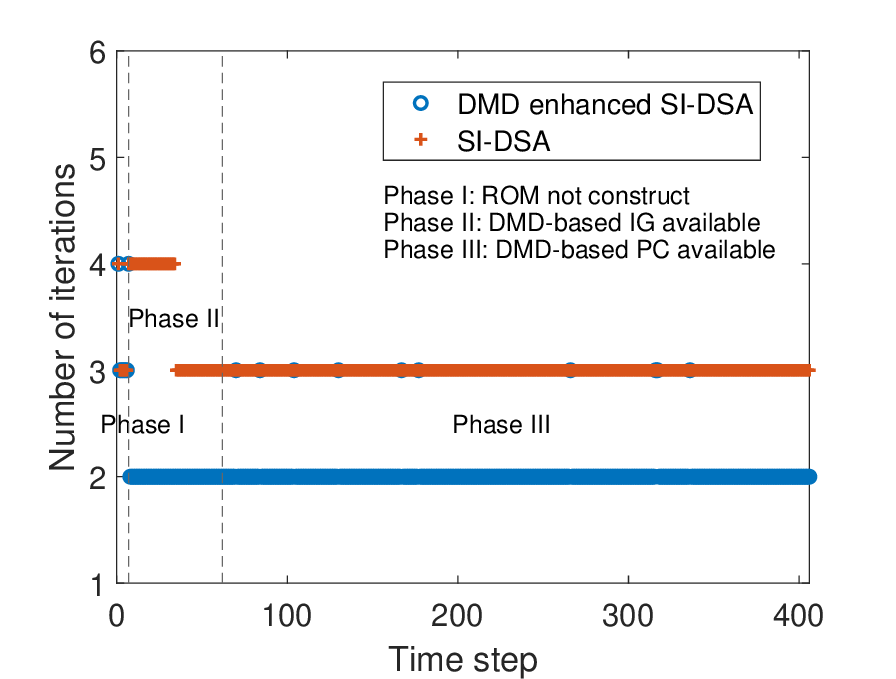}
        \end{minipage}
        \label{fig:scattering_rho}
    \end{subfigure}
     \hspace{-0.01\textwidth}
    \begin{subfigure}[b]{0.45\textwidth}
        \centering
        \hspace{-10mm}
        \begin{minipage}{6cm}
        \includegraphics[scale=0.55]{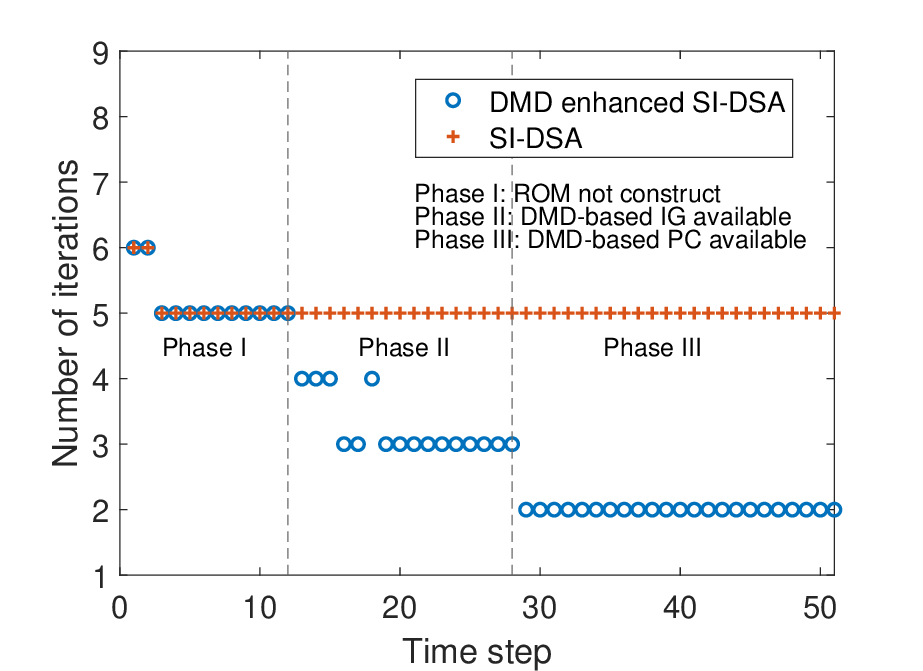}
        \end{minipage}
        \label{fig:scattering_rho_y0}
    \end{subfigure}
    \caption{Number of iterations in each time step with condition $\sigma_s=1$. Left: CFL=0.25. Right: CFL=2.}    \label{fig:Gaussian_source_num_iters}
\end{figure}

\textbf{Influence of time step size.} We first investigate the impact of the time step size with $\sigma_s=1$.  We set the time step size as  $\Delta t=\textrm{CFL}\Delta x$ with $\textrm{CFL}=0.25,0.5,1$ and $2$. The relative computational time for different time step size of each phases are  shown in Tab. \ref{tab:Gaussian_source_diffCFL}. The results of historical iteration counts with CFL=0.25 and CFL=2 are shown in Fig. \ref{fig:Gaussian_source_num_iters}. 
We find that the higher the time step size, $\Delta t$, is, the greater computational saving is gained by the DMD enhanced SI-SA in Phase III.  We also note that the time step size impacts the portion of each phases in time-marching, which further affects the overall speedup gained. 

\textbf{Influence of scattering strength.}
The results for $\sigma_s=0.1,1,10$ and $100$ are shown in Tab. \ref{tab:Gaussian_source_diffs}.  The overall speedup gained by our DMD-based acceleration is approximately $1.37\times$ for $\sigma_s=0.1$, $1.81\times$ for $\sigma_s=1$, $1.80\times$ for $\sigma_s=10$ and  $1.87\times$ for $\sigma_s=100$. In conclusion, our DMD enhanced SI-DSA gain greater computational saving in the intermediate ($\sigma_s=1,10$) and scattering dominant regime ($\sigma_s=100$).   

\begin{table}[htbp]
\centering
\caption{The average relative computational time for the isotropic scattering problem in Sec. \ref{sec:iso} with different time step size $\Delta t=\textrm{CFL}\Delta x$, where the DMD truncation tolerance of the initial guess is: $\epsilon_{\rm IG} = 10^{-9}$, the tolerance of preconditioner is: $\epsilon_{\rm PC} = 10^{-6}$, and the tolerance of adaptive update is: $\epsilon_{\rm UP} = 10^{-9}$. Phase II: ROM-enhanced initial guess available. Phase III: ROM-enhanced initial guess and preconditioners available.}

\begin{subtable}[t]{1.0\textwidth}
\centering
\begin{tabular}{ccccc}
\toprule
 &$\textrm{CFL}=0.25$& $\textrm{CFL}=0.5$ & $\textrm{CFL}=1$& $\textrm{CFL}=2$ \\
\midrule
Phase II &72.00\% & 77.20\% & 66.71\% &  66.72\% \\
Phase III  & 79.92\% & 64.15\% & 49.26\% & 47.11\%  \\
\bottomrule\\
\end{tabular}
\end{subtable}
\label{tab:Gaussian_source_diffCFL}
\end{table}

\begin{table}[htbp]
\centering
\caption{The average relative computational time for the isotropic scattering problem in Sec. \ref{sec:iso} with different $\sigma_s$. The DMD truncation tolerance of the initial guess is: $\epsilon_{\rm IG} = 10^{-9}$, the tolerance of preconditioner is: $\epsilon_{\rm PC} = 10^{-6}$, and the tolerance of adaptive update is: $\epsilon_{\rm UP} = 10^{-9}$. Phase II: ROM-enhanced initial guess available. Phase III: ROM-enhanced initial guess and preconditioners available.}
\begin{subtable}[t]{1.0\textwidth}
\centering
\scalebox{0.78}{\begin{tabular}{cccccccc}
\toprule
\multicolumn{1}{c}{}& \multicolumn{2}{c}{Phase II} & \multicolumn{2}{c}{Phase III} & \multicolumn{1}{c}{Entire period} & \multicolumn{1}{c}{DMD precomputation} & \multicolumn{1}{c}{IG construction}\\
\cmidrule(lr){2-3} \cmidrule(lr){4-5}
&\small{Time Steps} & \small{Avg. relative Time} & \small{Time Steps} & \small{Avg. relative Time} & \small{Avg. relative Time} & \small{Avg. relative Time} & \small{Avg. relative Time}\\
\midrule
$\sigma_s=0.1$ & 17 &91.14\% & 75 & 65.17\%  & 72.96\% & 2.1\% & $3.26\times 10^{-4}\%$\\
$\sigma_s=1$ & 17 &66.71\% & 75 & 49.26\% & 55.61\% & 1.56\% & $2.16\times 10^{-4}\%$\\
$\sigma_s=10$ & 16 &63.20\% & 77 & 52.74\% & 58.68\% & 1.36\% & $1.78\times 10^{-4}\%$\\
$\sigma_s=100$ & 14 & 46.32\% & 81 & 50.62\% & 53.51\% & 0.92\% & $9.78\times 10^{-5}\%$\\
\bottomrule\\
\end{tabular}}
\end{subtable}
\label{tab:Gaussian_source_diffs}
\end{table}

\subsection{Variable scattering problem\label{sec:variable-scattering}}
\begin{figure}[htbp]
    \centering
    \begin{subfigure}[b]{0.33\textwidth}
        \centering
        \hspace{-20mm}
        \begin{minipage}{4cm}
        \includegraphics[scale=0.39]{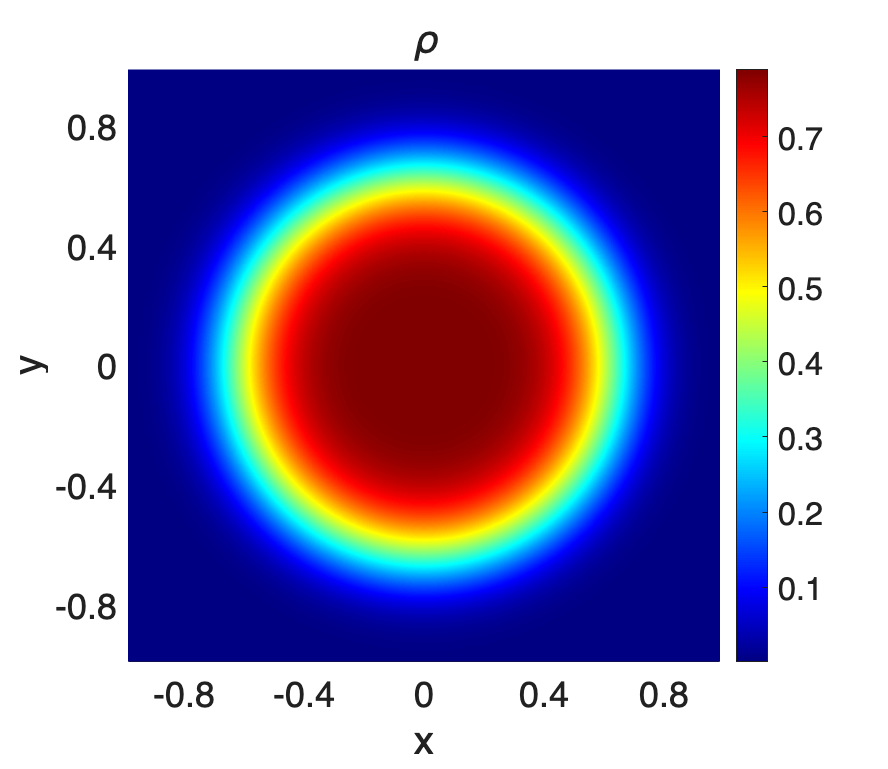}
        \end{minipage}
        \label{fig:scattering_rho}
    \end{subfigure}
     \hspace{-0.01\textwidth}
     \begin{subfigure}[b]{0.33\textwidth}
        \centering
        \hspace{-18mm}
        \begin{minipage}{4cm}
        \includegraphics[scale=0.39]{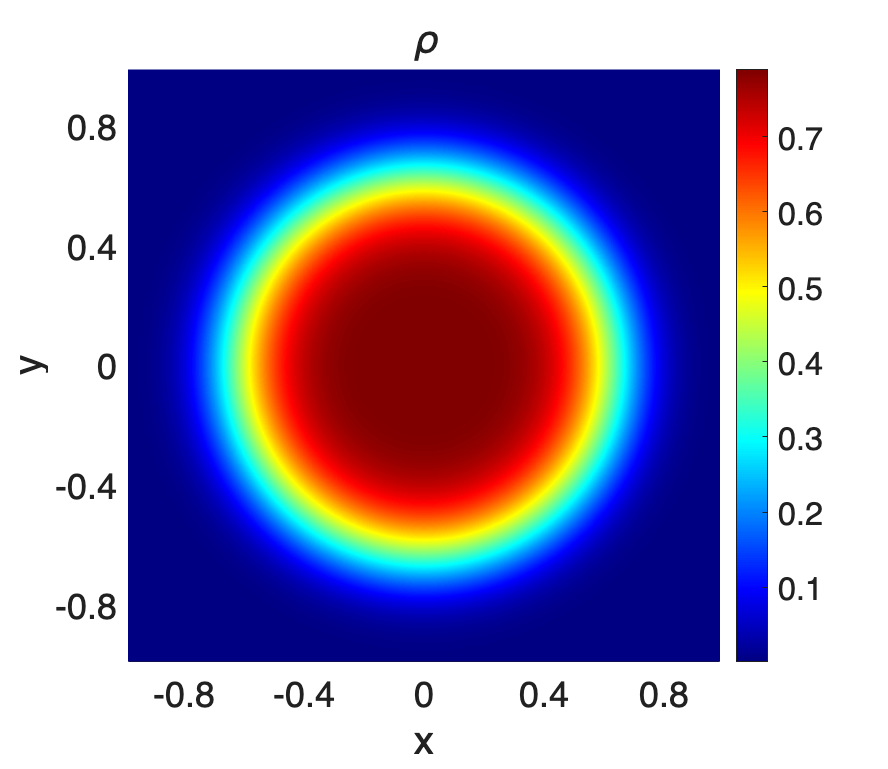}
        \end{minipage}
        \label{fig:scattering_rho}
    \end{subfigure}
     \hspace{-0.01\textwidth}
    \begin{subfigure}[b]{0.33\textwidth}
        \centering
        \hspace{-10mm}
        \begin{minipage}{4cm}
        \includegraphics[scale=0.34]{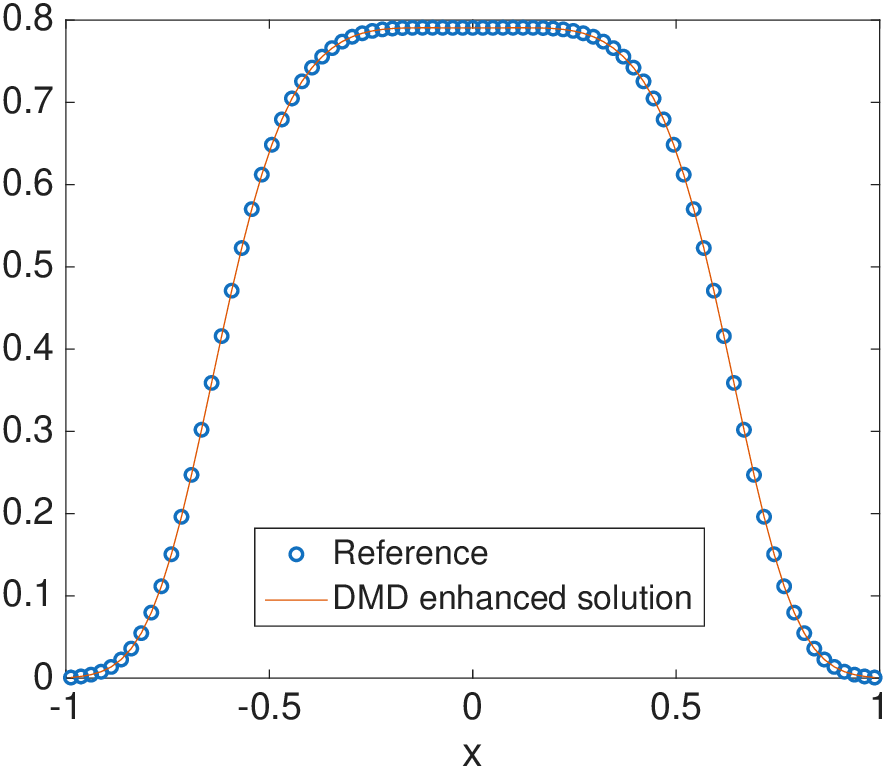}
        \end{minipage}
        \label{fig:scattering_rho_y0}
    \end{subfigure}
    \caption{ Density profile at $t=2.5$ for the variable scattering problem in Sec. \ref{sec:variable-scattering}. Left: solution by DMD enhanced SI-DSA. Middle: reference solution by SI-DSA. Right: the comparison of two solutions with $y=0$.
   }    \label{fig:lattice_density}
\end{figure}

We consider a 2D-2V problem on the computational domain $\Omega_{\bx}=[-1,1]^2$ with source $G(x,y)=0$, a variable scattering cross section, zero inflow boundary conditions and a smooth initial condition:
\begin{subequations}
    \begin{align}
        f(\bx, \bv,0) = \frac{1}{4\pi\zeta^2}{\rm exp}\Big(-\frac{x^2+y^2}{4\zeta^2}\Big), \;\;\zeta=10^{-2}, \;\; (x,y)\in[-1,1]^2,\\
        \sigma_s(x,y)=\left\{\begin{array}{ll}
    99.9c^4(c+\sqrt{2})^2(c-\sqrt{2})^2+0.1, & \text{if  }   c=\sqrt{x^2+y^2}<1, \\
    1, & \text{otherwise,}
    \end{array}\right.
    \end{align}
\end{subequations}
This test is a challenging multiscale problem: the scattering strength $\sigma_s(x,y)$ smoothly transits from $0.1$ to $100$ from the center of the region to the boundary. In other words, there is a smooth transition from the transport dominant regime to the scattering dominant regime.
We partition the computational domain with an $81\times81$ mesh, and set the time step size as $\Delta t=\Delta x$. 
We use CL($40,6$) quadrature in the angular space, hence  $N_{\bv}=240$. We solve the problem from $t=0$ to $t=2.5$, and present numerical solutions at  $t=2.5$  in Fig.~\ref{fig:lattice_density}. We observe that the solution given by SI-DSA and our ROM enhanced SI-SA matches well. In fact, the $l_2$ difference between these two final time numerical solution is approximately $1.30\times10^{-11}$.

\begin{figure}[htbp]
    \centering
    \begin{subfigure}[b]{0.49\textwidth}
        \centering
        \begin{minipage}{6cm}
        \hspace{-5mm}
        \includegraphics[scale=0.45]{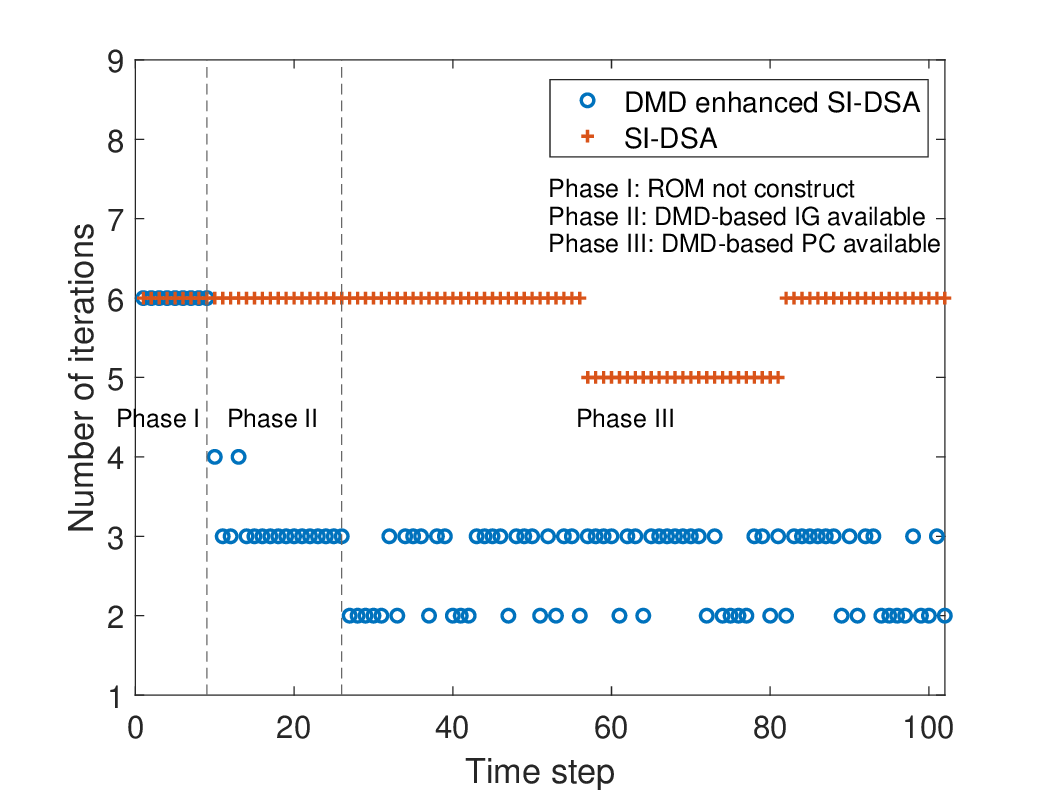}
        \end{minipage}
    \end{subfigure}
     \hspace{-0.01\textwidth}
    \begin{subfigure}[b]{0.49\textwidth}
        \centering
        \hspace{-15mm}
        \begin{minipage}{6cm}
        \includegraphics[scale=0.45]{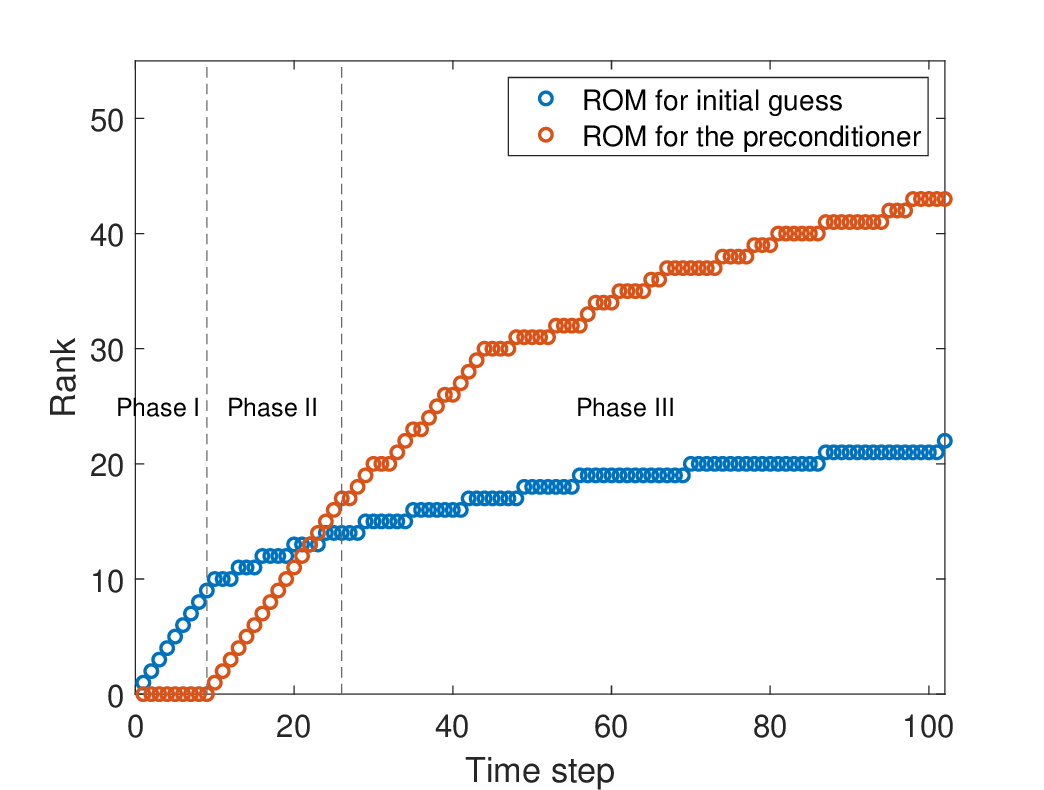}
        \end{minipage}
    \end{subfigure}
    \caption{Left: number of iterations in each time step for the variable scattering problem in Sec. \ref{sec:variable-scattering}. Right: rank in each time step for the variable scattering problem in Sec. \ref{sec:variable-scattering}.}    \label{fig:variable_scatt_num_iters}
\end{figure}

\textbf{Acceleration.}  
In Tab. \ref{tab:variable_scatter_time}, we summarize the number of time steps of various phases in our on-the-fly ROM-enhanced SI-SA, the average number of transport sweeps per time step and relative computational time in each phase or over the entire time marching. The computational time is reduced by approximately $40\%$ in phase II (ROM-enhanced initial guess available), $45\%$ in phase III (both ROM-enhanced initial guess and preconditioners avaible), and $40\%$ over the entire period. Fig. \ref{fig:variable_scatt_num_iters} shows that DMD enhanced SI-DSA requires significantly fewer iterations than the standard SI-DSA.

\textbf{Overhead in the construction of ROM and ROM-enhanced initial guesses.} Tab. \ref{tab:variable_scatter_DMDprepare} reports the relative overhead of constructing ROMs and ROM-enhanced initial guesses compared to the total time of implicit time marching with SI-DSA. We observe that the relative overhead of constructing ROMs and ROM-enhanced initial guesses  amounts to less than $2\%$, while resulting in substantial overall computational saving.

\textbf{Effective control of ROM rank.} The right panel of Fig. \ref{fig:variable_scatt_num_iters} shows the ROM rank for the initial guess and the preconditioner. As expected, the rank predictably rises linearly during construction, whereas it can be effectively controlled during the update stage via truncated SVD.

\begin{table}[htbp]
\centering
\caption{Average number of transport sweeps per time steps and average relative computational time for both methods in different phases for the variable scattering problem in Sec. \ref{sec:variable-scattering}. Phase I: both ROM-enhanced initial guesses and preconditioners are not available. Phase II: ROM-enhanced initial guess available. Phase III: ROM-enhanced initial guess and preconditioners available.}
\label{tab:variable_scatter_time}

\begin{subtable}[t]{1.0\textwidth}
\centering
\begin{tabular}{cccccc}
\toprule
& Time steps &SI-DSA & DMD enhanced SI-DSA & SI-DSA & DMD enhanced SI-DSA \\
\midrule
\multicolumn{2}{c}{} & \multicolumn{2}{c}{Avg.  sweeps per time step} & \multicolumn{2}{c}{Avg. relative computational time}\\
\cmidrule(lr){3-4} \cmidrule(lr){5-6}
Phase I  & 9  & 7 & Under construction & 100\% & Under construction \\ 
Phase II & 17 & 7 & 4.12  &  100\% &  60.38\% \\
Phase III & 76 & 6.67 & 3.58 &  100\% &  54.95\% \\
Over 
the entire period & 102 & 6.75 & 3.97 &  100\% &  60.35\% \\
\bottomrule
\end{tabular}
\end{subtable}
\end{table}

\begin{table}[h]
\centering
\caption{The total relative time compared to obtaining solutions for all time steps using SI-DSA for the variable scattering problem in Sec. \ref{sec:variable-scattering}. }
\begin{subtable}[t]{1.0\textwidth}
\centering
\begin{tabular}{cc}
\toprule
  Construction of ROMs &   Initial guess construction    \\
\midrule
1.29\% & $1.43\times 10^{-4}\%$   \\
\bottomrule\\
\end{tabular}
\end{subtable}
\label{tab:variable_scatter_DMDprepare}
\end{table}

\subsection{Lattice problem \label{sec:lattice}}
\begin{figure}[htbp]
    \centering
    \begin{subfigure}[b]{0.4\textwidth}
        \centering
        \begin{minipage}{6cm}
        \includegraphics[scale=0.51]{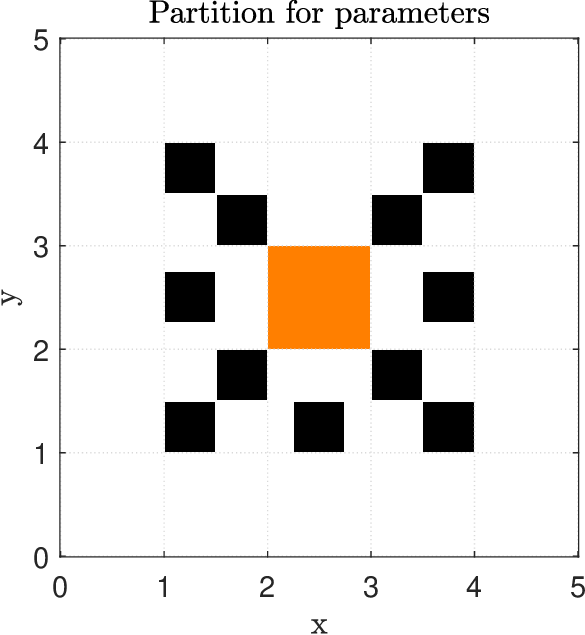}
        \end{minipage}
        \label{fig:lattice_parameter}
    \end{subfigure}
    \hspace{-0.04\textwidth}
    \begin{subfigure}[b]{0.45\textwidth}
        \centering
        \begin{minipage}{6cm}
        \includegraphics[scale=0.45]{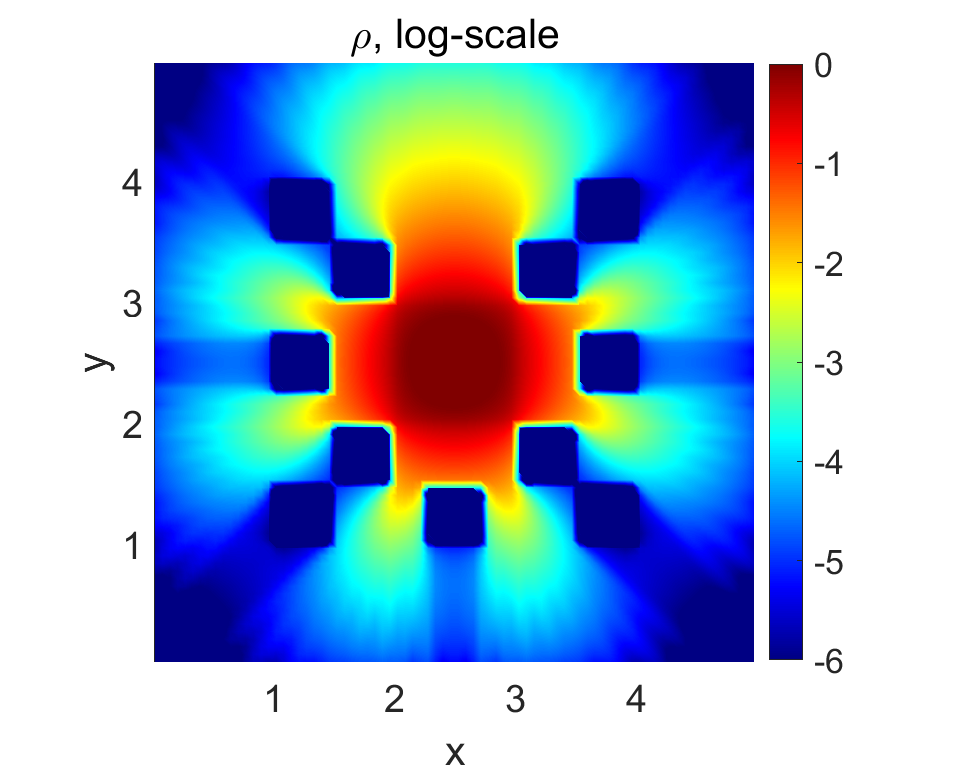}
        \end{minipage}
        \label{fig: lattice_log_solution}
    \end{subfigure}
    \caption{The set-up and a reference solution for the lattice problem in Sec. \ref{sec:lattice}. Left: the set-up of cross section for the lattice problem. Right: the DMD enhanced solution under log-scale.} 
    \label{fig: lattice_solution}
\end{figure}

The computational domain is a square $\mathcal{D}=[0,5]^{2}$ with zero inflow boundary conditions. The configuration of material properties is shown in the left picture of  Fig. \ref{fig: lattice_solution}. The black regions represent pure absorption zones with $(\sigma_{a},\sigma_{s}) = (100,0)$, while the remaining areas are pure scattering zones with $(\sigma_{a},\sigma_{s}) = (0,1)$.  A source term is imposed in the orange region, i.e.,
\begin{equation}
G(x, y)=\left\{\begin{array}{ll}
1.0, & \text{if } |x-2.5|<0.5 \text{ and } |y-2.5|<0.5, \\
0, & \text{otherwise.}
\end{array}\right.
\end{equation}

The initial condition is $f(\bx, \bv,0)=0$. We use a $80\times80$ uniform spatial grid to partition the computational domain, and apply CL($40,6$) quadrature rule for the angular discretization. The time step size is chosen as $\Delta t = \frac{1}{16}$, and we solve the problem from $t=0$ to $t=5$. We plot the reference density profile under the log-scale in  Fig. \ref{fig: lattice_solution}.

\begin{figure}[htbp]
    \centering
    \begin{subfigure}[b]{0.49\textwidth}
        \centering
        \begin{minipage}{6cm}
        \hspace{-5mm}
        \includegraphics[scale=0.45]{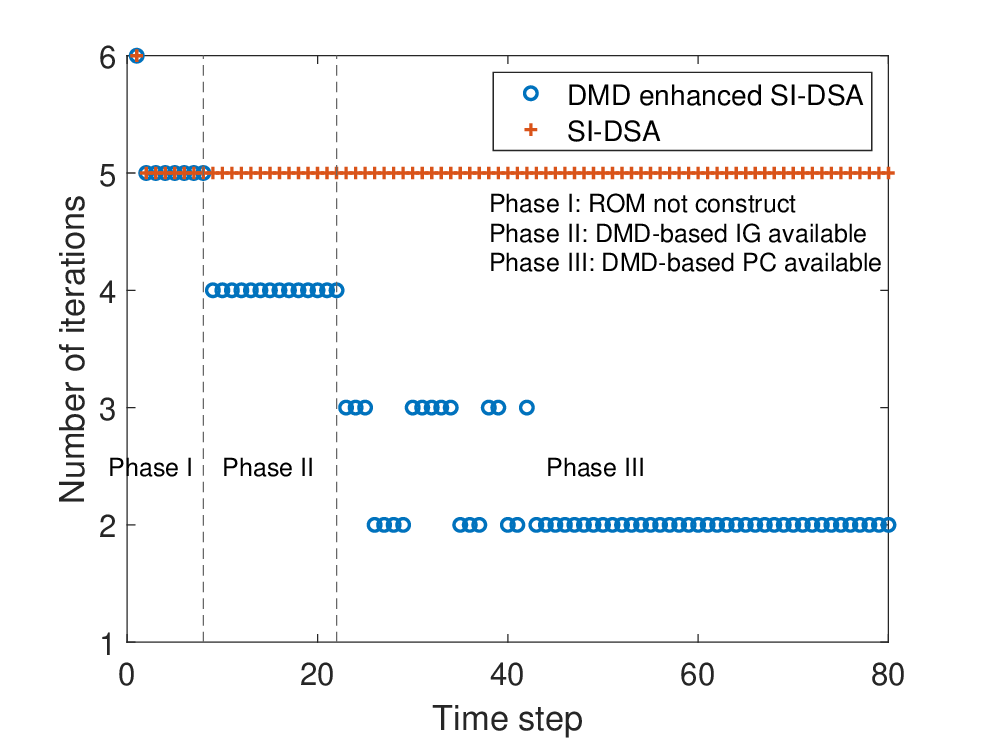}
        \end{minipage}
        \label{fig:scattering_rho}
    \end{subfigure}
     \hspace{-0.01\textwidth}
    \begin{subfigure}[b]{0.49\textwidth}
        \centering
        \hspace{-15mm}
        \begin{minipage}{6cm}
        \includegraphics[scale=0.45]{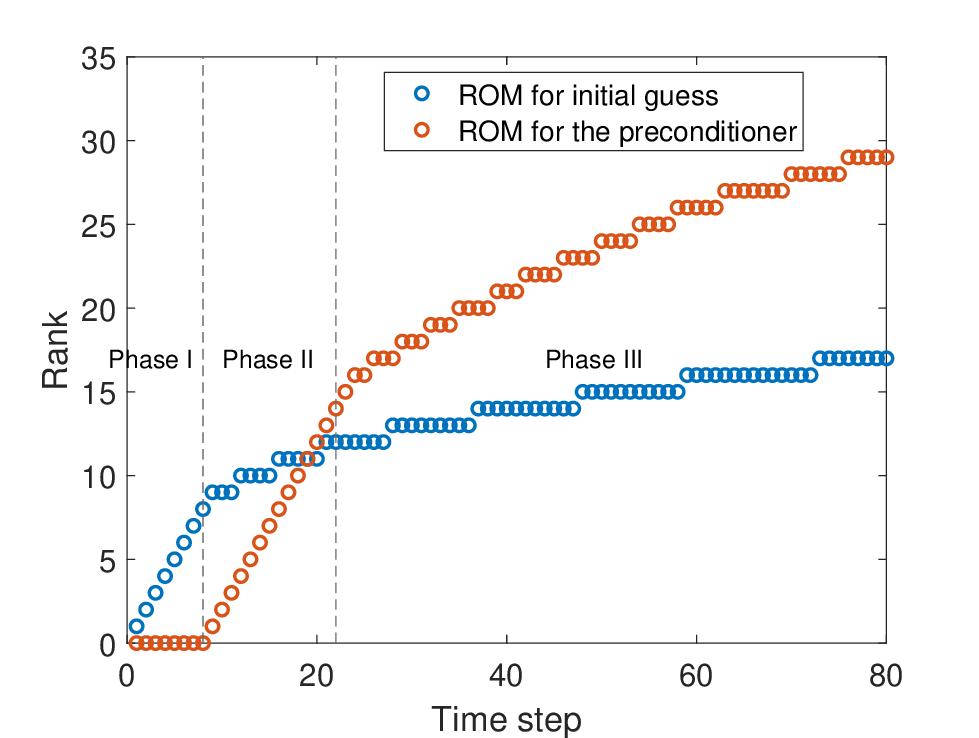}
        \end{minipage}
        \label{fig:scattering_rho_y0}
    \end{subfigure}
    \caption{Left: number of iterations in each time step for the lattice problem in Sec. \ref{sec:lattice}. Right: rank in each time step for the lattice problem in Sec. \ref{sec:lattice}.}    \label{fig:lattice_num_iters}
\end{figure}

\begin{enumerate}
\item \textbf{Acceleration.} In Tab. \ref{tab:lattice_cost} and Tab. \ref{tab:lattice_time}, we present the relative computational cost for $(\epsilon_{\textrm{IG}},\epsilon_{\textrm{PC}},\epsilon_{\textrm{UP}})=(10^{-6},10^{-6},10^{-6})$. Fig. \ref{fig:lattice_num_iters} further demonstrates that DMD enhanced SI-DSA can significantly reduce iteration counts compared to SI-DSA. The number of iterations can be reduced to $2$ iterations, resulting $1.82\times$ speedup in phase $3$ and an average $1.52\times$ acceleration over the entire $80$ steps time marching.

\item \textbf{Overhead for the construction of ROMs and initial guesses.} As presented
in  Tab. \ref{tab:lattice_time}, we observe that the relative time for building the DMD model and constructing initial guesses remains under $2\%$. In other words, with marginal overhead, great acceleration can be achieved. In addition, the overhead of constructing initial guesses is negligible compared to the time of building ROMs.

\item \textbf{Effective control of ROM rank.} The rank of the ROM for the initial guess and for the preconditioner is shown in the right panel of Fig.  \ref{fig:lattice_num_iters}. We observe that the rank linearly increases during the construction stage as expected, while it can be effectively controlled during the update stage using truncated SVD.
\end{enumerate}

\begin{table}[h]
\centering
\caption{Relative computational cost of SI-DSA and DMD enhanced SI-DSA for the lattice problem in Sec. \ref{sec:lattice}. Phase I: ROM not construct. Phase II: ROM-enhanced initial guess available. Phase III: ROM-enhanced initial guess and preconditioners available.}
\label{tab:lattice_cost}

\begin{subtable}[t]{1.0\textwidth}
\centering
\scalebox{0.9}{
\begin{tabular}{cccccc}
\toprule
\multicolumn{6}{c}{$\epsilon_{\mathrm{IG}} = 10^{-6}$, $\epsilon_{\mathrm{PC}} = 10^{-6}$, $\epsilon_{\mathrm{UP}} = 10^{-6}$}\\
\midrule
\multicolumn{2}{c}{}& \multicolumn{2}{c}{Avg. number of transport sweeps} & \multicolumn{2}{c}{Avg. realtive computational time}\\
\cmidrule(lr){3-4} \cmidrule(lr){5-6}
& Time steps &SI-DSA & DMD enhanced SI-DSA & SI-DSA & DMD enhanced SI-DSA \\
Phase I  & 8  & 6.13 & Under construction & 100\% & Under construction \\ 
Phase II & 14 & 6 & 5 &  100\% & 86.51\% \\
Phase III & 58 & 6 & 3.19 &  100\% & 54.94\% \\
Over 
the entire period & 80 & 6.01 & 3.8 &  100\% &  65.47\% \\
\bottomrule
\end{tabular}}
\end{subtable}
\end{table}

\begin{table}[h]
\centering
\caption{The total relative time compared to obtaining solutions for all time steps using SI-DSA for the lattice problem in Sec. \ref{sec:lattice}. }
\label{tab:lattice_time}
\begin{subtable}[t]{1.0\textwidth}
\centering
\begin{tabular}{cc}
\toprule
\multicolumn{2}{c}{$\epsilon_{\mathrm{IG}} = 10^{-6}$, $\epsilon_{\mathrm{PC}} = 10^{-6}$, $\epsilon_{\mathrm{UP}} = 10^{-6}$}\\
\midrule
\multicolumn{1}{c}{Construction of ROMs} & \multicolumn{1}{c}{Initial guess construction}\\
1.71\% &  $1.63\times10^{-4}\%$\\
\bottomrule
\end{tabular}
\end{subtable}
\end{table}

\subsection{Comparison and combination with the DMD-based preconditioner in \cite{mcclarren2022data} \label{sec:Gaussian_combine}}

In this subsection, we compare the proposed method with the DMD preconditioner in \cite{mcclarren2022data} and combine our approach with this preconditioner by replacing DSA with it, namely ``Hybrid with McClarren-Haut method". This test demonstrates the flexibility of the proposed acceleration strategy in being integrated with preconditioners beyond DSA.

The numerical example considered in this section is the same as that in Sec. \ref{sec:iso}, the simulation time remains from $t=0$ to $t=2.5$, and the parameters are set as: $\sigma_s=1$, CFL=1, $\epsilon_{\textrm{SISA}}=10^{-8}$, $\epsilon_{\textrm{IG}}=10^{-6}$, $\epsilon_{\textrm{PC}}=10^{-6}$, and $\epsilon_{\textrm{UP}}=10^{-6}$.

\begin{figure}[htbp]
    \centering
    \begin{subfigure}[b]{0.45\textwidth}
        \centering
        \hspace{-10mm}
        \begin{minipage}{6cm}
        \includegraphics[scale=0.45]{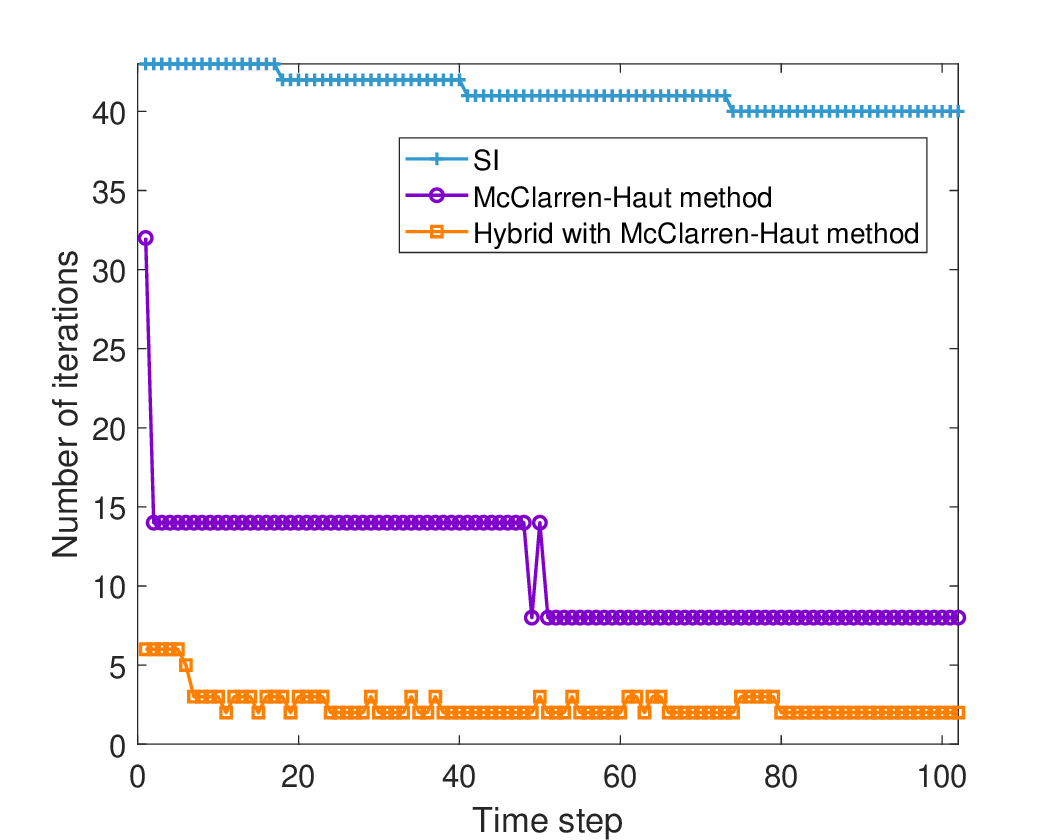}
        \end{minipage}
        \label{fig:scattering_rho}
    \end{subfigure}
    \begin{subfigure}[b]{0.45\textwidth}
        \centering
        \begin{minipage}{6cm}
        \includegraphics[scale=0.45]{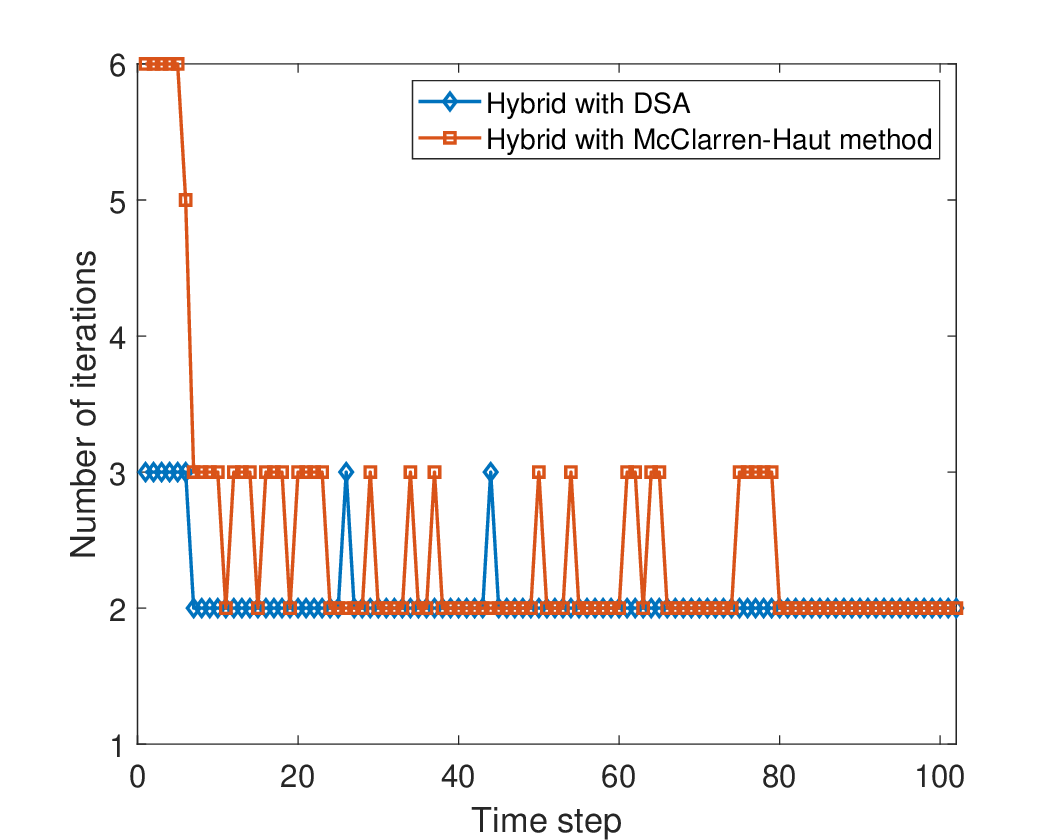}
        \end{minipage}
        \label{fig:scattering_rho_y0}
    \end{subfigure}
    \caption{Number of iteration steps of different methods in Sec. \ref{sec:Gaussian_combine}. Left: comparison of SI, McClarren-Haut method, and Hybrid with McClarren-Haut method. Right: comparison of Hybrid with DSA and Hybrid with McClarren-Haut method.}  \label{fig:Gaussian_source_comparison}
\end{figure}

Combining our approach and the McClarren-Haut method, enjoys the advantages of the two: our proposed method provides a more accurate initial guess and the first iteration, and McClarren-Haut method provides a surrogate model for acceleration in the following iteration steps (see Fig. \ref{fig:Gaussian_source_comparison}). In the right panel of Fig. \ref{fig:Gaussian_source_comparison}, it can be observed that enhancing McClarren-Haut method with our approach achieves performance comparable to enhancing SI-DSA with our approach. Therefore, this hybrid approach can serve as an efficient acceleration strategy in scenarios where conventional DSA is not applicable, e.g. target problems in \cite{mcclarren2022data}.

\begin{figure}[htbp]
    \centering
    \begin{subfigure}[b]{0.45\textwidth}
        \centering
        \hspace{-10mm}
        \begin{minipage}{6cm}
        \includegraphics[scale=0.45]{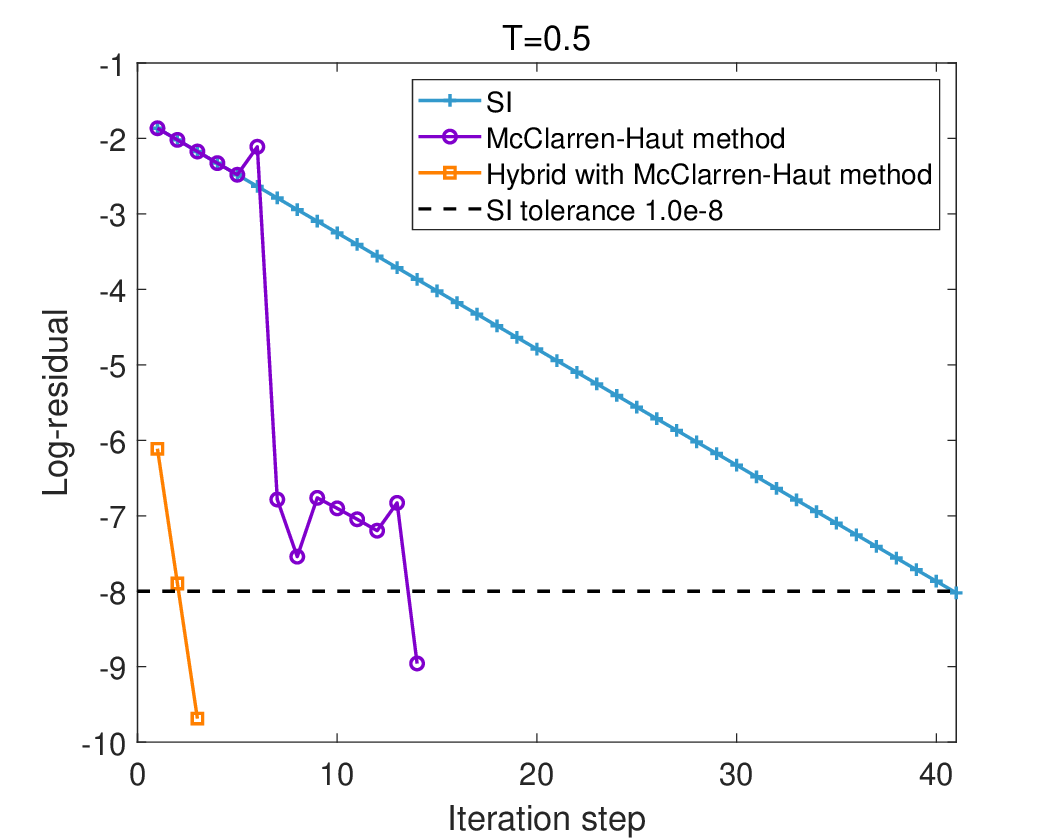}
        \end{minipage}
        \label{fig:scattering_rho}
    \end{subfigure}
    \begin{subfigure}[b]{0.45\textwidth}
        \centering
        \begin{minipage}{6cm}
        \includegraphics[scale=0.45]{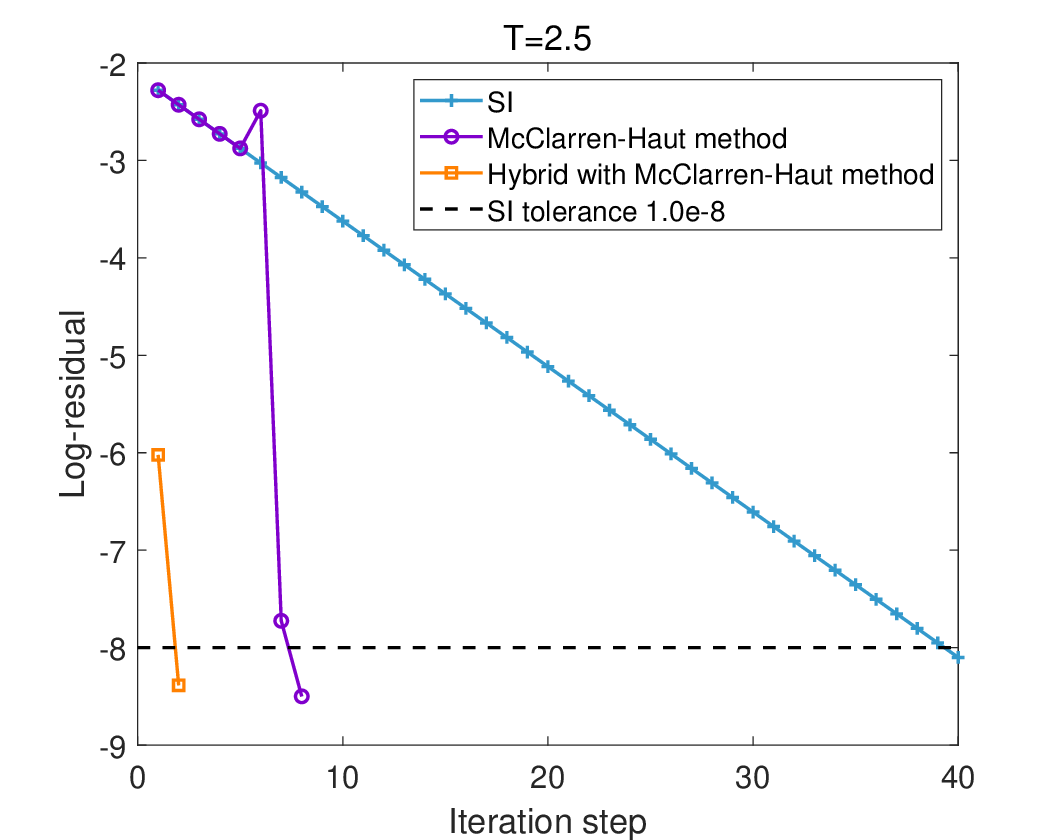}
        \end{minipage}
        \label{fig:scattering_rho_y0}
    \end{subfigure}
    \caption{ Log-residual history for the test in Sec. \ref{sec:Gaussian_combine}. Left: T=0.5. Right: T=2.5. }    \label{fig:Gaussian_source_residual}
\end{figure}

Fig. \ref{fig:Gaussian_source_residual} displays the log-residual history at $t=0.5$ and $t=2.5$. Oscillations exists in the residual history of McClarren-Haut method. Combining it with our acceleration strategy,  the residuals declines more rapidly and does not oscillate.

\section{Conclusion}
\label{sec:conclusion}
Leveraging low-rank structures across time and the kinetic description of RTE, we propose an offline-free, on-the-fly ROM-based acceleration of SI-DSA for the implicit time marching of  RTE. 

Our method consists of three phases: (1) construction of the ROM to provide enhanced initial guesses, (2) construction of a second ROM to enhance DSA preconditioner, and (3) application and adaptive updates of ROM-enhanced initial guesses and preconditioners. To accommodate the streaming nature of solution data in implicit time marching, memory-efficient reformulation of RTE, DMD, incremental SVD, and cheap error indicators are leveraged to construct and update ROMs in an efficient and memory-lean manner.

In our numerical tests, the proposed method effectively reduces the   number of iterations for convergence by approximately $40\%$ compared to  standard SI-DSA and achieves approximately $1.4\times$ to $2\times$ speedup for all examples, while paying marginal overhead to construct and update ROMs. Furthermore, we demonstrate the flexibility of our method by integrating it with preconditioners other than DSA. 

In the future, we plan to extend the proposed approach to nonlinear thermal radiations, and multi-group problem. 

\section*{Acknowledgement}
This work was partially supported by the Hong Kong Research Grants Council grants Early Career Scheme 26302724 and General Research Fund 16306825.

\section*{CRediT authorship contribution statement}

{\bf Ningxin Liu:} Writing – original draft, Writing – review \& editing, Visualization, Validation, Software, Methodology, Data curation, Conceptualization.

{\bf Zhichao Peng:} Writing – review \& editing, Visualization, Validation, Software, Methodology, Data curation, Conceptualization, Supervision, Funding acquisition.

\section*{Declaration of generative AI and AI-assisted technologies in the writing process}
During the preparation of this work, the authors used AI to check for grammar errors and improve readability. After using this tool/service, the authors reviewed and edited the content as needed and take full responsibility for the content of the publication.

\appendix
\section{Chebyshev-Legendre quadrature rule\label{appx:cl}}
The Chebyshev-Legendre (CL) quadrature for integrals on $\mathbb{S}^2$ can be seen as the tensor product of the Chebyshev rule for the unit cycle and the Gaussian-Legendre rule for the z-component of the angular direction $\bv_z\in[-1,1]$. The  quadrature points and weights of the normalized  $N_{\phi}$-points Chebshev quadrature rule for the unit circle is 
\begin{equation*}
\left\{(\phi_j,\omega_{\phi,j}): \;\phi_j = \frac{(2j-1)\pi}{N_\phi}, \;\;\omega_{\phi,j}=\frac{1}{N_\phi},\; j=1,\dots,N_\phi\right\}.
\end{equation*}
We denote the quadrature points and the quadrature weights of the normalized $N_{v_z}$-points Gauss-Legendre rule as $\{(v_{z,j},\omega_{z,j})\}_{j=1}^{N_{v_z}}$. The quadrature points and weights of the CL($N_{\phi}$,$N_{v_z}$) quadrature rule are defined as:
\begin{equation*}
\bv_{j}=\left(\cos(\phi_{j_1})\sqrt{1-v_{z,j_2}^2},\;\;\sin(\phi_{j_1})\sqrt{1-v_{z,j_2}^2},\;\;v_{z,j_2}\right), \quad\;\omega_j=\omega_{\phi,j_1}\omega_{z,j_2},
\end{equation*}
where  $1\leq j_1\leq N_\phi$, $1\leq j_2\leq N_{v_z}$, and $1\leq j=(j_2-1) N_{\phi}+j_1\leq N_{\bv}$ with $N_{\bv}=N_{v_z}N_{\phi}$.
\section{Incremental SVD\label{appx:i-svd}}
Here, we briefly review the low-rank incremental SVD update algorithm in \cite{brand2006fast}.
Given $\BR_{m-1}\approx \BU\BS\BV^\rT$, where $\BU\in\mathbb{R}^{N_{\bx}\times r},\BS\in\mathbb{R}^{r\times r},\BV\in\mathbb{R}^{m\times r}$. Define $\bq=\frac{1}{||(\BI-\BU\BU^T)\brho||_2}(\BI-\BU\BU^T)\brho$, then $||\bq||_2=1$ and $\BU^T\bq=\boldsymbol{0}$.  
\begin{align}
   \BR_m=\begin{bmatrix}\BR_{m-1} & \brho^m\end{bmatrix} 
   &\approx \begin{bmatrix}\BU\BS\BV^\rT&\brho^m\end{bmatrix} 
   =\begin{bmatrix}\BU\BS\BV^\rT&\BU\BU^T\brho^m+\bq\end{bmatrix}\notag\\
   &= \begin{bmatrix}\BU&\bq\end{bmatrix}
   \begin{bmatrix}
   \BS\BV^\rT &  \BU^\rT\brho\\
    \mathbf{0}       &   ||\BI-\BU\BU^T\rho||_2
   \end{bmatrix}=\underbrace{\begin{bmatrix}\BU&\bq\end{bmatrix}}_{N_{\bx}\times(r+1)}
   \underbrace{\begin{bmatrix}
   \BS &  \BU^\rT\brho\\
    \mathbf{0}       &   ||\BI-\BU\BU^T\rho||_2
   \end{bmatrix}}_{(r+1)\times(r+1)}
   \underbrace{\begin{bmatrix}
   \BV^\rT       & \mathbf{0}\\
   \mathbf{0}    & 1 
   \end{bmatrix}}_{(r+1)\times m}.
   \label{eq:svd_core}
\end{align}
The full or truncated SVD of the core $(r+1)\times(r+1)$ matrix can be computed as
\begin{equation}
\begin{bmatrix}
   \BS &  \BU^\rT\brho\\
    \mathbf{0}       &   ||\BI-\BU\BU^T\rho||_2
   \end{bmatrix}
   \approx\BU_{\textrm{core}}\BS_{\textrm{core}}\BV_{\textrm{core}}^T
\end{equation}
with $O(r^3)$ cost.
Hence, truncated SVD of $\BR_m$ can be approximated as
\begin{equation}
\BR_m\approx\hat{\BU}\hat{\BS}\hat{\BV}^\rT,\;\text{where}\;\hat{\BU}=\begin{bmatrix}\BU & \bq\end{bmatrix}\BU_{\textrm{core}},\;
\hat{\BS}=\BS_{\textrm{core}},\;
\hat{\BV}=\begin{bmatrix}
\BV & \mathbf{0}\\
\mathbf{0} & 1
\end{bmatrix}\BV_{\textrm{core}}.
\end{equation}
The total cost of computing the truncated SVD of $\BR_m$ is $O((N_{\bx}+m+r)r^2)$. We note that necessary truncation may be applied in equation \eqref{eq:svd_core}. 

If one compute the truncated SVD of a $N_{\bx}\times m$ matrix iteratively using the above procedure with $O((N_{\bx}+m+r)r^2)$ in each step, the total cost is $O(N_{\bx}mr)$.

\bibliographystyle{elsarticle-num} 
\bibliography{ref}

@article{choi2021space,
  title={Space--time reduced order model for large-scale linear dynamical systems with application to Boltzmann transport problems},
  author={Choi, Youngsoo and Brown, Peter and Arrighi, William and Anderson, Robert and Huynh, Kevin},
  journal={Journal of Computational Physics},
  volume={424},
  pages={109845},
  year={2021},
  publisher={Elsevier}
}

@article{buchan2015pod,
  title={A POD reduced order model for resolving angular direction in neutron/photon transport problems},
  author={Buchan, Andrew G and Calloo, AA and Goffin, Mark G and Dargaville, Steven and Fang, Fangxin and Pain, Christopher C and Navon, Ionel Michael},
  journal={Journal of Computational Physics},
  volume={296},
  pages={138--157},
  year={2015},
  publisher={Elsevier}
}

@article{peng2022reduced,
  title={A reduced basis method for radiative transfer equation},
  author={Peng, Zhichao and Chen, Yanlai and Cheng, Yingda and Li, Fengyan},
  journal={Journal of Scientific Computing},
  volume={91},
  number={1},
  pages={5},
  year={2022},
  publisher={Springer}
}

@article{peng2024micro,
  title={A micro-macro decomposed reduced basis method for the time-dependent radiative transfer equation},
  author={Peng, Zhichao and Chen, Yanlai and Cheng, Yingda and Li, Fengyan},
  journal={Multiscale Modeling \& Simulation},
  volume={22},
  number={1},
  pages={639--666},
  year={2024},
  publisher={SIAM}
}

@article{tano2021affine,
  title={Affine reduced-order model for radiation transport problems in cylindrical coordinates},
  author={Tano, Mauricio and Ragusa, Jean and Caron, Dominic and Behne, Patrick},
  journal={Annals of Nuclear Energy},
  volume={158},
  pages={108214},
  year={2021},
  publisher={Elsevier}
}

@article{behne2022minimally,
  title={Minimally-invasive parametric model-order reduction for sweep-based radiation transport},
  author={Behne, Patrick and Vermaak, Jan and Ragusa, Jean C},
  journal={Journal of Computational Physics},
  volume={469},
  pages={111525},
  year={2022},
  publisher={Elsevier}
}

@article{mcclarren2019calculating,
  title={Calculating time eigenvalues of the neutron transport equation with dynamic mode decomposition},
  author={McClarren, Ryan G},
  journal={Nuclear Science and Engineering},
  volume={193},
  number={8},
  pages={854--867},
  year={2019},
  publisher={Taylor \& Francis}
}

@article{coale2023reduced,
  title={Reduced order models for thermal radiative transfer problems based on moment equations and data-driven approximations of the Eddington tensor},
  author={Coale, Joseph M and Anistratov, Dmitriy Y},
  journal={Journal of Quantitative Spectroscopy and Radiative Transfer},
  volume={296},
  pages={108458},
  year={2023},
  publisher={Elsevier}
}

@article{chen2021low,
  title={A low-rank Schwarz method for radiative transfer equation with heterogeneous scattering coefficient},
  author={Chen, Ke and Li, Qin and Lu, Jianfeng and Wright, Stephen J},
  journal={Multiscale Modeling \& Simulation},
  volume={19},
  number={2},
  pages={775--801},
  year={2021},
  publisher={SIAM}
}

@article{kopp1963synthetic,
  title={Synthetic method solution of the transport equation},
  author={Kopp, HJ},
  journal={Nuclear Science and Engineering},
  volume={17},
  number={1},
  pages={65--74},
  year={1963},
  publisher={Taylor \& Francis}
}

@article{adams1992diffusion,
  title={Diffusion synthetic acceleration of discontinuous finite element transport iterations},
  author={Adams, Marvin L and Martin, William R},
  journal={Nuclear science and engineering},
  volume={111},
  number={2},
  pages={145--167},
  year={1992},
  publisher={Taylor \& Francis}
}

@article{wareing1993new,
  title={New diffusion-sythetic acceleration methods for the S [sub N] equations with corner balance spatial differencing},
  author={Wareing, TA},
  year={1993}
}

@article{peng2024reduced,
  title={Reduced order model enhanced source iteration with synthetic acceleration for parametric radiative transfer equation},
  author={Peng, Zhichao},
  journal={Journal of Computational Physics},
  volume={517},
  pages={113303},
  year={2024},
  publisher={Elsevier}
}

@article{mcclarren2022data,
  title={Data-driven acceleration of thermal radiation transfer calculations with the dynamic mode decomposition and a sequential singular value decomposition},
  author={McClarren, Ryan G and Haut, Terry S},
  journal={Journal of Computational Physics},
  volume={448},
  pages={110756},
  year={2022},
  publisher={Elsevier}
}

@article{schmid2010dynamic,
  title={Dynamic mode decomposition of numerical and experimental data},
  author={Schmid, Peter J},
  journal={Journal of fluid mechanics},
  volume={656},
  pages={5--28},
  year={2010},
  publisher={Cambridge University Press}
}

@book{pomraning2005equations,
  title={The equations of radiation hydrodynamics},
  author={Pomraning, Gerald C},
  year={2005},
  publisher={Courier Corporation}
}

@article{adams2001discontinuous,
  title={Discontinuous finite element transport solutions in thick diffusive problems},
  author={Adams, Marvin L},
  journal={Nuclear science and engineering},
  volume={137},
  number={3},
  pages={298--333},
  year={2001},
  publisher={Taylor \& Francis}
}

@article{guermond2010asymptotic,
  title={Asymptotic analysis of upwind discontinuous Galerkin approximation of the radiative transport equation in the diffusive limit},
  author={Guermond, Jean-Luc and Kanschat, Guido},
  journal={SIAM Journal on Numerical Analysis},
  volume={48},
  number={1},
  pages={53--78},
  year={2010},
  publisher={SIAM}
}

@article{adams2002fast,
  title={Fast iterative methods for discrete-ordinates particle transport calculations},
  author={Adams, Marvin L and Larsen, Edward W},
  journal={Progress in nuclear energy},
  volume={40},
  number={1},
  pages={3--159},
  year={2002},
  publisher={Elsevier}
}

@article{alcouffe1977diffusion,
  title={Diffusion synthetic acceleration methods for the diamond-differenced discrete-ordinates equations},
  author={Alcouffe, Raymond E},
  journal={Nuclear Science and Engineering},
  volume={64},
  number={2},
  pages={344--355},
  year={1977},
  publisher={Taylor \& Francis}
}

@article{lorence1989s,
  title={An s 2 synthetic acceleration scheme for the one-dimensional sn equations with linear discontinuous spatial differencing},
  author={Lorence Jr, Leonard J and Morel, JE and Larsen, Edward W},
  journal={Nuclear Science and Engineering},
  volume={101},
  number={4},
  pages={341--351},
  year={1989},
  publisher={Taylor \& Francis}
}

@article{ramone1997transport,
  title={A transport synthetic acceleration method for transport iterations},
  author={Ramone, Gilles L and Adams, Marvin L and Nowak, Paul F},
  journal={Nuclear science and engineering},
  volume={125},
  number={3},
  pages={257--283},
  year={1997},
  publisher={Taylor \& Francis}
}

@article{gol1964quasi,
  title={A quasi-diffusion method of solving the kinetic equation},
  author={Gol'Din, V Ya},
  journal={USSR Computational Mathematics and Mathematical Physics},
  volume={4},
  number={6},
  pages={136--149},
  year={1964},
  publisher={Elsevier}
}

@article{anistratov1993nonlinear,
  title={Nonlinear methods for solving particle transport problems},
  author={Anistratov, Dmitri Yu and Gol'Din, Vladimir Ya},
  journal={Transport Theory and Statistical Physics},
  volume={22},
  number={2-3},
  pages={125--163},
  year={1993},
  publisher={Taylor \& Francis}
}

@article{olivier2023family,
  title={A family of independent variable Eddington factor methods with efficient preconditioned iterative solvers},
  author={Olivier, Samuel and Pazner, Will and Haut, Terry S and Yee, Ben C},
  journal={Journal of Computational Physics},
  volume={473},
  pages={111747},
  year={2023},
  publisher={Elsevier}
}

@article{schmid2011applications,
  title={Applications of the dynamic mode decomposition},
  author={Schmid, Peter J and Li, Larry and Juniper, Matthew P and Pust, Oliver},
  journal={Theoretical and computational fluid dynamics},
  volume={25},
  number={1},
  pages={249--259},
  year={2011},
  publisher={Springer}
}

@phdthesis{tu2013dynamic,
  title={Dynamic mode decomposition: Theory and applications},
  author={Tu, Jonathan H},
  year={2013},
  school={Princeton University}
}

@article{schmid2011application,
  title={Application of the dynamic mode decomposition to experimental data},
  author={Schmid, Peter J},
  journal={Experiments in fluids},
  volume={50},
  number={4},
  pages={1123--1130},
  year={2011},
  publisher={Springer}
}

@article{brand2006fast,
  title={Fast low-rank modifications of the thin singular value decomposition},
  author={Brand, Matthew},
  journal={Linear algebra and its applications},
  volume={415},
  number={1},
  pages={20--30},
  year={2006},
  publisher={Elsevier}
}

@article{tang2025synthetic,
  title={Synthetic Acceleration Preconditioners for Parametric Radiative Transfer Equations based on Trajectory-Aware Reduced Order Models},
  author={Tang, Ning and Peng, Zhichao},
  journal={arXiv preprint arXiv:2509.05001},
  year={2025}
}

@article{ren2019fast,
  title={{A fast algorithm for radiative transport in isotropic media}},
  author={Ren, Kui and Zhang, Rongting and Zhong, Yimin},
  journal={Journal of Computational Physics},
  volume={399},
  pages={108958},
  year={2019},
  publisher={Elsevier}
}

@article{peng2025flexible,
  title={{A flexible GMRES solver with reduced order model enhanced synthetic acceleration preconditioner for parametric radiative transfer equation}},
  author={Peng, Zhichao},
  journal={Journal of Computational Physics},
  volume={534},
  pages={114004},
  year={2025},
  publisher={Elsevier}
}

@article{smith2023variable,
  title={{Variable dynamic mode decomposition for estimating time eigenvalues in nuclear systems}},
  author={Smith, Ethan and Variansyah, Ilham and McClarren, Ryan},
  journal={Nuclear Science and Engineering},
  volume={197},
  number={8},
  pages={1769--1778},
  year={2023},
  publisher={Taylor \& Francis}
}

@article{chen2008ifem,
  title={{iFEM: an innovative finite element methods package in MATLAB}},
  author={Chen, Long},
  journal={Preprint, University of Maryland},
  volume={20},
  year={2008}
}

@article{larsen2009advances,
  title={{Advances in discrete-ordinates methodology}},
  author={Larsen, Edward W and Morel, Jim E},
  journal={Nuclear computational science: A century in review},
  pages={1--84},
  year={2009},
  publisher={Springer}
}

@article{matsuda2025reduced,
  title={{Reduced basis methods for parametric steady-state radiative transfer equation}},
  author={Matsuda, Kimberly and Chen, Yanlai and Cheng, Yingda and Li, Fengyan},
  journal={Journal of Computational Physics},
  pages={114597},
  year={2025},
  publisher={Elsevier}
}

@article{tencer2017accelerated,
  title={{Accelerated solution of discrete ordinates approximation to the boltzmann transport equation for a gray absorbing--emitting medium via model reduction}},
  author={Tencer, John and Carlberg, Kevin and Larsen, Marvin and Hogan, Roy},
  journal={Journal of Heat Transfer},
  volume={139},
  number={12},
  pages={122701},
  year={2017},
  publisher={American Society of Mechanical Engineers}
}

\end{document}